\documentclass[11pt]{amsart}
\usepackage{amsmath,amscd,amssymb}
\usepackage[mathscr]{eucal}
\usepackage{amsmath}
\usepackage{enumerate}
\usepackage{color}
\usepackage{mathtools}

\DeclareMathOperator{\arccot}{arccot}

 \newtheorem{thm}{Theorem}[section]
 \newtheorem{cor}[thm]{Corollary}
 \newtheorem{lemma}[thm]{Lemma}
 \newtheorem{prop}[thm]{Proposition}

 \theoremstyle{definition}

  \newtheorem{exas}[thm]{Examples}
	
 \theoremstyle{remark}
 \newtheorem{rem}[thm]{Remark}

\newtheorem{remark}[thm]{Remark}

\numberwithin{equation}{section}

\def\R{\mathbb{R}}
\def\C{\mathbb{C}}
\def\N{\mathbb{N}}
\def\Z{\mathbb{Z}}
\def\a{\alpha}
\def\b{\beta}
\def\ep{\epsilon}
\def\l{\lambda}

\def\Re#1{\operatorname{Re}#1}
\def\Sect#1{\operatorname{Sect}#1}
\def\Rat{\mathcal{R}_0(\C_+)}

\def\L{\L}

\def\B{\mathcal{B}}
\def\Bes{\mathcal{B}}
\def\Bov{\mathcal{E}}
\def\Bq{{\Bes_0}}
\def\A{\mathcal{A}}
\def\D{\mathcal{D}}
\def\E{\mathcal{E}}
\def\F{\mathcal{F}}
\def\fti{\mathcal{F}^{-1}}
\def\G{\mathcal{G}}
\def\L{\mathcal{L}}
\def\LT{\mathcal{LM}}

\def\lt{\mathcal{L}}
\def\V{\mathcal{V}}
\def\W{\mathcal{W}}

\def\RR{\mathcal{R}}
\def\H{\mathcal{H}}

\def\Hol{\operatorname{Hol}}

\begin{document}

\title[Bounded functional calculi for unbounded operators]{Bounded functional calculi for unbounded operators}

\author{Charles Batty}
\address{St. John's College\\
University of Oxford\\
Oxford OX1 3JP, UK
}

\email{charles.batty@sjc.ox.ac.uk}

\author{Alexander Gomilko}
\address{Faculty of Mathematics and Computer Science\\
Nicolas Copernicus University\\
Chopin Street 12/18\\
87-100 Toru\'n, Poland
}

\email{alex@gomilko.com}

\author{Yuri Tomilov}
\address{
Institute of Mathematics\\
Polish Academy of Sciences\\
\'Sniadeckich 8\\
00-656 Warsaw, Poland
}

\email{ytomilov@impan.pl}

\begin{abstract}
This article summarises the theory of several bounded functional calculi for unbounded operators that have recently been discovered.   The extend the Hille--Phillips calculus for (negative) generators $A$ of certain bounded $C_0$-semigroups, in particular for bounded semigroups on Hilbert spaces and bounded holomorphic semigroups on Banach spaces.  They include  functions outside the Hille-Phillips class, and they generally give sharper bounds for the norms of the resulting operators $f(A)$.    The calculi are mostly based on appropriate reproducing formulas for the relevant classes of functions, and they rely on significant and interesting developments of function theory.   They are compatible with standard functional calculi and they admit appropriate convergence lemmas and spectral mapping theorems.   They can also be used to derive several well-known operator norm-estimates, provide generalisations of some of them, and extend the general theory of operator semigroups.   Our aim is to help readers to make use of these calculi without having to understand the details of their construction.
\end{abstract}




\date \today

\maketitle

\section{Introduction}

Functional calculus for operators $A$ has a long history, originating with self-adjoint operators on Hilbert spaces, considered by Hilbert himself in the case of bounded operators,  and by von Neumann for unbounded operators.   In that context, $f(A)$ can be defined as a bounded operator for every bounded measurable function $f$ on the spectrum of $A$. 
For other operators on Hilbert spaces or Banach spaces, a functional calculus is a way of converting one operator (often unbounded) into a coherent collection of associated bounded operators.  Many applications may require $f(A)$ to be defined for only a few functions $f(z)$, such as $\exp(-tz)$ for operator semigroups $(e^{-tA})_{t\ge0}$,  and $(z-1)^n(z+1)^{-n}$ for powers of cogenerators.  
Results concerning specific functions can be discussed individually as has often been the case, but constructing functional calculi covering many functions has advantages.   
A unified approach is likely to be efficient in the long run, and the treatment of a functional calculus is likely to be more attractive mathematically than considering many different functions individually.  
 Moreover, including more functions $f$ often improves the sharpness of estimates of the norm of $f(A)$, and may possibly widen the scope of applications in the future.

For operators which are not self-adjoint, a functional calculus usually has to be restricted to holomorphic functions $f$ defined on an open set $\Omega$ containing the spectrum $\sigma(A)$ of $A$.  If $A$ is bounded, the Riesz-Dunford calculus can be defined by
\[
f(A) = \frac{1}{2\pi i} \int_\gamma f(z) (z-A)^{-1} \, dz,
\]
where $\gamma$ is a contour within $\Omega$ around $\sigma(A)$.   Extending this to unbounded operators runs into complications, as the contours would naturally have to go to infinity.   Sometimes the contour $\gamma$ may have to pass through a point in the boundary of $\Omega$,  and then convergence of the integral above may depend on further assumptions on $f$ and $A$. 

If $A$ is an injective sectorial operator of angle $\theta \in [0,\pi)$, one can circumvent this problem by using the extended sectorial calculus as described in the book of Haase \cite{HaaseB}.  The prototype for functional calculi of this type was a paper by Bade on strip-type operators \cite{Bade}.   If $f \in H^\infty(\Sigma_\psi)$ for some $\psi \in (\theta,\pi)$, $f(A)$ can be defined as a closed operator on $X$.  Many differential operators $A$ have a bounded $H^\infty$-calculus, in the sense that $f(A)$ is bounded for all $f \in H^\infty(\Sigma_\psi)$, but there are operators which do not have bounded $H^\infty$-calculus (see \cite[Remark 4.7]{BGTad}, for example).    If $-A$ is the generator of a bounded $C_0$-semigroup on $X$, then $A$ is sectorial of angle $\pi/2$ and $\sigma(A) \subseteq \overline\C_+$.  In this case it may be more natural to consider a half-plane calculus as in \cite{BHM}.  

For the sectorial (or half-plane) calculus, the functions $f$ have to be defined on a larger sector than the sectorial angle of $A$ (or a larger half-plane than $\C_+$), and this is sometimes an unnatural restriction.  If $A$ is sectorial of angle strictly less than $\pi/2$ (written as $A \in \Sect(\pi/2-))$, we may consider functions $f \in H^\infty(\C_+)$.  Then $f(A)$ is defined in the sectorial calculus, but not necessarily bounded.    Vitse \cite{V1} showed that $f(A)$ is a  bounded operator for all $f$ in the Banach algebra $\Bes$ of ``analytic Besov functions'' on $\C_+$ which is continuously embedded in $H^\infty(\C_+)$ (see Section \ref{B-calc} and Theorem 4.1).   Thus every operator $A \in \Sect(\pi/2-)$ has a bounded $\B$-calculus.   Subsequently Schwenninger \cite{Sch1} improved and extended Vitse's results, including an estimate of the norms of operators of the form $f(A)e^{-tA}$ for small $t>0$ and $f \in H^\infty(\C_+)$.  


An early example of a bounded functional calculus for an unbounded operator $A$ is the Hille-Phillips calculus.   This calculus defines $f(A)$ whenever $f$ is the Laplace transform of a bounded Borel measure on $\R_+$  and $-A$ is the generator of a bounded $C_0$-semigroup on a Banach space.  It was described in \cite[Chapter XV]{HP} and it played a substantial role in developing semigroup theory.   

The use of analytic Besov functions for functional calculus goes back to Peller \cite{Pel} in the discrete case.  He showed that if $T$ is a power-bounded operator on a Hilbert space, then $\|p(T)\| \le C\|p\|_{\B_d}$ for all polynomials $p(z)$, where $\B_d$ is the analogue of $\B$ for the unit disc.  Since the polynomials are dense in $\B_d$, it follows that $T$ has a bounded functional calculus for $\B_d$.   In the continuous-parameter case, the analogous class of analytic Besov functions on $\C_+$ is more complicated (for example, it is not separable), and it has received little attention in the literature.   The earliest research which successfully considered functional calculus for negative generators of bounded $C_0$-semigroups on Hilbert space was carried out by White in his doctoral thesis \cite{White}, but it remained virtually unknown for many years.    He proved estimates of the form $\|f(A)\| \le C\|f\|_\B$ for functions $f$ in a subclass of the Hille-Phillips algebra.    There was a further contribution by Haase \cite{Haase} who (inspired by \cite{V1}, and unaware of White's work) obtained a similar estimate for a larger class of functions. 

Many other functional calcuii have been proposed and developed by numerous contributors.   This article is not a comprehensive survey of the subject.  Instead, we set out the basic properties of a small group of other functional calculi which have recently been developed and are related to the calculus in \cite{V1}.  A bounded $\Bes$-calculus for a class of operators $A$, including the negative generators of all bounded $C_0$-semigroups on Hilbert spaces and all bounded holomorphic $C_0$-semigroups, is described in Section 3, based on our papers \cite{BGT} and \cite{BGT2}.  
Section 4 gives an account of two classes of function spaces, $\D_s$ and $\H_\psi$, and the associated $\D$-calculus and $\H$-calculus for operators $A\in\Sect(\pi/2-)$, based on \cite{BGT3}.   Spaces in both these classes contain $\B$, and they also included functions which are unbounded on $\C_+$ or only defined on proper subsectors of $\C_+$.
Another bounded calculus for a Banach algebra $\A$ of functions which can be applied negative generators of bounded $C_0$-semigroups on Hilbert spaces has been constructed by Arnold and Le Merdy \cite{ALM} and it is briefly described in Section \ref{A-calc}.     

The aim of this article is to present the main results about these calculi in a coherent fashion which will enable other researchers to appreciate their value and to apply them when possible, without having to piece together information from several lengthy papers.   We describe the definitions and main properties of the algebras and the calculi, how they relate to each other and to the sectorial and half-plane calculi, and how they relate to some problems.  Proofs are omitted or briefly sketched. In general terms, the proofs in our papers depend mainly on complex analysis and functional analysis, with some aspects of harmonic analysis.  


We present the functional calculi only for negative generators of bounded $C_0$-semigroups, but they can be adapted to unbounded semigroups.   Assume that $A$ is the negative generator of a $C_0$-semigroup $(T(t))_{t\ge0}$ and $\|T(t)\| \le Me^{\omega t}$, and  $X$ is a Hilbert space or $(T(t))_{t\ge0}$ is a holomorphic semigroup.   Then $f(A)$ can be defined to be $g(A-\omega)$, where $g(z) := f(z+\omega)$, provided that $g$ belongs to $\B$, $\D_s$ or $\H_\psi$, as appropriate.
   
In the IWOTA conference, the first author of this article lectured mostly on the $\B$-calculus, with only a brief mention of the $\D$- and $\H$-calculi.   The first author of \cite{ALM} lectured on the $\A$-calculus.
 
\subsection* {Notation}  
We will use the following notation:
\begin{enumerate}[\null\hskip1pt]
\item $\R_+ :=[0,\infty)$,
\item $\C_+ := \{z \in\C: \Re z>0\}$, $\overline{\C}_+ = \{z \in\C: \Re z\ge0\}$,
\item $\Sigma_\theta := \{z\in\C: z \ne 0, |\arg z|<\theta\}$ for $\theta \in (0,\pi]$, $\Sigma_0 := (0,\infty)$,
\item $L(X,Y)$ is the space of bounded linear operators from one Banach space $X$ to another $Y$; we write $L(X)$ when $Y=X$.
\end{enumerate}

\section{Semigroup generators and Hille-Phillips calculus} \label{HP-calc}
In this paper $A$ will be a closed operator on a Banach space $X$ with a dense domain $D(A)$, and its spectrum $\sigma(A)$ will be contained in $\overline\Sigma_\theta$ for some $\theta\in[0,\pi)$.   Except for material on the $\H$-calculus, $\theta$ will be in $[0,\pi/2]$, and $-A$ will generate a bounded $C_0$-semigroup $(T(t))_{t\ge0}$.  Then we define
\[
K_A := \sup_{t\ge0} \|T(t)\|.
\]

Let $\mathfrak A$ be a Banach algebra which is continuously included in $H^\infty(\Sigma_\theta)$ where $\theta \in (0,\pi]$, and assume that $\mathfrak A$ contains the functions $r_z(\l):= (\l+z)^{-1}$ for $\l\in\Sigma_\theta, \, z\in \Sigma_{\pi-\theta}$.   An {\it $\mathfrak A$-calculus} for $A$ is a bounded algebra homomorphism $\Phi : \mathfrak A \to L(X)$ such that $\Phi(r_z) = (z+A)^{-1}$ for all $z\in\Sigma_{\pi-\theta}$.

An operator $A$ on a Banach space $X$ is {\it sectorial of angle} $\theta \in [0,\pi)$ if $\sigma(A) \subset \overline \Sigma_\theta$ and, for each $\psi \in (\theta,\pi)$,
\begin{equation} \label{R01}
M_\psi(A) := \sup_{z\in\Sigma_{\pi-\psi}} \|z(z+A)^{-1}\| < \infty.
\end{equation}

We write $\Sect(\theta)$ for the class of all sectorial operators of angle $\theta$ for $\theta \in [0,\pi)$ on Banach spaces, and $\Sect(\pi/2-) := \bigcup_{\theta \in [0,\pi/2)} \Sect(\theta)$.   Then $A \in \Sect(\pi/2-)$ if and only if  $-A$ generates a (sectorially) bounded holomorphic $C_0$-semigroup on $X$.  
We refer to \cite{HaaseB} for the general theory of sectorial operators, and to \cite[Section 3.7]{ABHN} for the theory of holomorphic semigroups.

Let $A$ be an operator and  assume that $\sigma(A) \subset \overline \C_+$ and
\begin{equation}\label{R11}
M_A := M_{\pi/2}(A) = \sup_{z\in \C_+} \|z(z+A)^{-1}\| < \infty.
\end{equation}
It follows from Neumann series that $\sigma(A) \subset \Sigma_\theta \cup \{0\}$ and $M_\theta(A) \le 2M_A$,
where
\[
\theta:=\arccos(1/(2M_A))<\pi/2.
\]
So $A \in \Sect(\theta) \subset \Sect(\pi/2-)$.   Conversely, if $A \in \Sect(\theta)$ where $\theta \in [0,\pi/2)$, then \eqref{R11} holds, with $M_A$ equal to $M_{\pi/2}(A)$ in (\ref{R01}).   Thus $-A$ generates a bounded holomorphic semigroup if and only if $\sigma(A) \subset \overline\C_+$ and (\ref{R11}) holds.  In that case, $M_A$ is a basic quantity associated with $A$, known as the {\it sectoriality constant} of $A$.

Let $M(\R_+)$ be the space of all bounded Borel measures on $\R_+$.  Then $M(\R_+)$ is a Banach algebra under convolution of measures.   For $\mu \in M(\R_+)$, let $\L\mu$ be the Laplace transform of $\mu$, and let 
\[
\LT = \{\L\mu : \mu \in M(\R_+)\}.
\]
This is a subalgebra of $H^\infty(\C_+)$,  and the map $\mu \mapsto \L \mu$ is an injective algebra homomorphism from $M(\R_+)$ into $H^\infty(\C_+)$.   Thus $\LT$ becomes a Banach algebra in the norm $\|\L\mu\|_{\text{HP}} := \|\mu\|_{M(\R_+)}$.

The following functions are in $\LT$:
\[
e_t(z):= e^{-tz}, \;t\ge0; \quad r_\l(z) := (\l+z)^{-1},\; \l\in\C_+;  \quad \frac{z-1}{z+1};  \quad \frac{1-e^{-z}}{z}.
\]

Assume that $-A$ is the generator of a bounded $C_0$-semigroup $(T(t))_{t\ge0}$ on a Banach space $X$.   For $f = \L\mu \in \LT$, define 
\begin{equation} \label{HP}
\Pi_A(f)x = \int_0^\infty T(t)x \,d\mu(t), \qquad x \in X.
\end{equation}
Then $\Pi_A$ is an $\LT$-functional calculus for $A$, with $T(t) = \Pi_A(e_t) = \exp(-tA)$ and $\|\Pi_A\| \le K_A$.   Moreover it is the unique  $\LT$-calculus for $A$, and it is known as the {\it Hille-Phillips calculus} or HP-{\it calculus}.

The Hille-Phillips calculus is compatible with the sectorial functional calculus in the sense that $\Pi_A(f) = f(A)$ whenever $f(A)$ is defined in the sectorial calculus \cite[Section 3.3]{HaaseB}.    Moreover there is a spectral inclusion theorem \cite[Theorem 2.7.4]{HaaseB}, and a convergence lemma for the HP-calculus \cite[Proposition 5.1.4]{HaaseB}.

We now discuss two questions which can be addressed via the HP-calculus.  Both questions are considered in \cite{Sasha} and we will return to them in later sections of this article. 

\subsection*{Cayley transform question}
Assume that $-A$ generates a bounded $C_0$-semigroup on a Banach space $X$ and let $V(A)$ be the cogenerator, so $V(A) = (A-1)(A+1)^{-1}$.  Let $V$ be the Cayley transform, so $V(z) = (z-1)(z+1)^{-1}$.   Then $V$ is the Laplace transform of $\delta_0 - 2 e^{-t}dt \in M(\R_+)$, so $V \in \LT$ and $\|V\|_{\text{HP}} = 3$.   The spectral radius of $V(A)$ is at most 1, and this raised a question about sharp upper bounds for the sequence  $(\|V(A)^n\|)_{n\in\N}$ for large $n$.     For the general case, it was shown in \cite{Thomee} that $\|V^n\|_{\text{HP}} = O(n^{1/2})$, so the Hille-Phillips calculus establishes that $\|V(A)^n\| = O(n^{1/2})$.   This estimate is known to be optimal for arbitrary $A$ \cite[Section 9]{Sasha}, so $\|V^n\|_{\text{HP}} \asymp n^{1/2}$.  We shall return to the Cayley transform question for operators on Hilbert space in Section 3 (Application 1)  and in Section 5.    

\subsection*{Inverse generator question}
Assume that $-A$ generates a bounded $C_0$-semigroup on $X$, and $A$ has dense range in $X$.  Then $A$ is injective and the operator $A^{-1}$ is densely defined.  The question was raised in \cite{deL} whether $-A^{-1}$ also generates a $C_0$-semigroup.  If $A \in \Sect(\pi/2-)$ then $A^{-1} \in \Sect(\pi/2-)$, so the answer is positive in this case.   For $t>0$ and $z \in \overline{\C}_+ \setminus \{0\}$, let $f_t(z) = e^{-t/z}$.   If operators $f_t(A)$ can be defined as bounded operators in some functional calculus for $A$, then the question should have a positive answer.   However $f_t$ cannot be defined continuously at $0$, so $f_t$ is not in $\LT$.  
We shall return to this question in Section \ref{B-calc} (Application 2).

\section{$\Bes$-calculus}  \label{B-calc}

To motivate the definition of the $\B$-calculus, we start by reconsidering the definition of the HP-calculus in Section \ref{HP-calc}.  Let $-A$ be the generator of a bounded $C_0$-semigroup $(T(t))_{t\ge0}$.  Let $f = \L\mu$, where $\mu\in M(\R_+)$ and $\mu(\{0\})=0$.   For $x \in X$, we may argue formally that
\begin{align*}
f(A)x &= \int_0^\infty T(t)x \, d\mu(t) \\
&= 4 \int_0^\infty \int_0^\infty \a t^2 e^{-2\a t} T(t)x \,d\a\,d\mu(t) \\
&= \frac{2}{\pi} \int_0^\infty \int_0^\infty \a t e^{-\a t} \int_\R (\a - i\b +A)^{-2}x \, e^{-i\b t} \, d\b \,d\a\,d\mu(t) \\
&= \frac{2}{\pi} \int_0^\infty \a \int_\R (\a - i\b +A)^{-2}x \int_0^\infty t e^{-(\a+i\b)t} \,d\mu(t) \,d\b\,d\a\\
&= -\frac{2}{\pi}  \int_0^\infty \a  \int_\R (\a - i\b +A)^{-2}x \, f'(\a+i\b) \,d\b\,d\a.
\end{align*}
Here we have assumed that the Fourier inversion formula can be applied to the function $t \mapsto \pi t e^{-\a t} T(t)x$ on $\R_+$ (extended to $\R$ by $0$) and its Fourier transform $\b \mapsto (\a +i\b + A)^{-2}x$, and that Fubini's theorem is valid.   The first assumption can be weakened by working in the weak operator topology.   Then both assumptions can be satisfied by making assumptions on $f$ and $A$ which involve two spaces $\B$ and $\E$ of holomorphic functions on $\C_+$.

\subsection*{The space $\Bes$}

Let $\Bes$ be the space of all holomorphic functions on $\C_+$ such that
\begin{equation} \label{Bnorm}
\|f\|_{\Bq} := \int_0^\infty \sup_{\b\in\R} |f'(\a+i\b)| \, d\a < \infty.
\end{equation}

The literature on this space $\Bes$ was very limited before the appearance of \cite{V1} and \cite{BGT},  and the corresponding space of  functions on the unit disc had received rather more attention.  The spaces are sometimes described as spaces of analytic Besov functions.    The origin for this terminology will be explained in Section 5 of this article.    We now set out some of the properties of individual functions in $\Bes$ (see \cite[Proposition 2.2]{BGT}).

\begin{prop}  \label{besprop}
Let $f \in \Bes$. 
\begin{enumerate}[\rm1.]
\item $f(\infty):= \lim_{\Re{z}\to\infty} f(z)$ exists in $\C$.
\item $f$ is bounded, and $\|f\|_\infty \le |f(\infty)| + \|f\|_\Bq$.
\item $f(i\b) := \lim_{\a\to0+} f(\a+i\b)$ exists, uniformly for $\b \in \R$.
\item The extended function $f$ is uniformly continuous on $\overline{\C}_+$.
\end{enumerate}
\end{prop}

We let $\Bq := \{f\in\B : f(\infty)=0\}$, so $\B = \B_0 \oplus \C$, where $\C$ represents the constant functions.

A norm on $\Bes$ is defined by
\[
\|f\|_\Bes := \|f\|_\infty + \|f\|_{\Bq}.
\]
The space $\Bes$ equipped with this norm has the following properties  (see \cite[Proposition 2.3, Lemma 2.14]{BGT}).

\begin{prop} \label{balg}
\begin{enumerate}[\rm1.]
\item The closed unit ball $U$ of $\Bes$ is compact in the topology of uniform convergence on compact subsets of $\C_+$.
\item  $\Bes$ is a Banach algebra, and $\Bes$ is not separable.
\item  The norm $\|\cdot\|_\Bq$ is equivalent to $\|\cdot\|_\B$ on $\Bq$, and  $\Bq$ is a closed ideal in the Banach algebra $\B$.
\item  $\LT \subset \Bes$ with continuous inclusion;  $\LT$ is not closed in $\Bes$.  
\item  $\LT$ is not dense in $\Bes$ in the norm-topology, but $\LT$ is dense in $\Bes$ in the topology of uniform convergence on compact subsets of $\C_+$.
\end{enumerate}
\end{prop}

\begin{exas} \label{exB}
1. The function $e^{-1/z}$ is not in $\Bes$, as it is not uniformly continuous near $0$.   However $z(z+1)^{-1}e^{-1/z}$ is in $\LT$ (see the proof of \cite[Theorem 3.3]{deL2}).   

\noindent
2. The function $\arccot z$ is defined on $\C_+$ by
\begin{equation} \label{arccot}
\arccot z = \frac{1}{2i} \log \left(\frac{z+i}{z-i}\right), \qquad z\in\C_+.
\end{equation}
This function is not bounded on $\C_+$, so it is not in $\Bes$.   The function $\exp(\arccot z)$ is bounded on $\C_+$, but it is not uniformly continuous near $i$, so $\exp(\arccot z)$ is not in $\Bes$ \cite[Example 3.19]{BGT3}.  

\noindent
3.   The following two functions are in $\B$ but not in $\LT$:
\[
e^{-z}(1+z)^i, \qquad e^{-2z} \int_0^1 \frac{\sinh(zt)}{t} \, dt.
\]
\end{exas}

Other classes of functions which are contained in $\Bes$ may be found in \cite[Section 3]{BGT}.  We now give a description of all functions in $\B$.

Let $\Hol(\C_+)$ be the space of all holomorphic functions on $\C_+$, and $\mathcal W$ be the space of all (equivalence classes of) measurable functions $g$ on $\C_+$ such that
\begin{equation}\label{A}
\|g\|_{\mathcal W}:=\int_0^\infty \underset{\beta\in \mathbb R}{\mathrm{ess\,sup}}\,|g(\alpha+i\beta)|\,d\alpha<\infty,
\end{equation}
with the norm given by (\ref{A}).    Then $\mathcal W$ is a Banach space, and
\begin{equation}\label{wb0}
\mathcal{W} \cap \operatorname{Hol}(\C_+) = \{ f' : f \in \Bes\}, \qquad \|f'\|_{\mathcal W} = \|f\|_\Bq.
\end{equation}

For $g \in \mathcal{W}$, let
 \begin{align} \label{defQ}
(Q g)(z):&= - \frac{2}{\pi}\int_0^\infty \alpha\int_\R \frac{g(\alpha+i\beta)}{(\alpha-i\beta+z)^2}\,d\beta\,d\alpha.
\end{align}
This defines an operator $Q$ which serves to obtain a reproducing formula for functions in $\B$ and to construct functions in $\B$ from arbitrary functions in $\W$, as in the following proposition.
The first statement is proved in \cite[Proposition 2.20]{BGT}, and the second and third statements are proved in \cite[Proposition 3.1]{BGT2}.

\begin{prop} \label{BHP2}
\begin{enumerate}[\rm1.]
\item  
Let $f \in \Bes$.   Then $f'\in\W$ and
\begin{equation} \label{Brepf}
f(z) = f(\infty) + Q(f')(z) = f(\infty) - \frac{2}{\pi} \int_0^\infty \a \int_\R \frac{f'(\a+i\b)}{(\a-i\b+z)^2}  \,d\b\,d\a.
\end{equation}
\item Let $g \in \W$.   Then $Qg \in \B_0$.  
\item  The operator $Q \in L(\W,\Bes_0)$, and it maps $\W \cap \Hol(\C_+)$ bijectively onto $\Bq$.
\end{enumerate}
\end{prop}

The Banach algebra $\Bes$ satisfies the following Convergence Lemma \cite[Lemma 8.1]{BGT2}.

\begin{lemma}\label{lem9.1}
Let $(f_k)_{k\ge 1}\subset \mathcal{B}$ be such that
$
\sup_{k\ge 1}\,\|f_k\|_{\mathcal{B}}<\infty.
$
Assume that for every $z\in \C_{+}$ there exists
\begin{equation*}\label{Conv21}
f(z):=\lim_{k\to\infty}\,f_k(z) \in \C,
\end{equation*}
and for every $r>0$, 
\begin{equation*}\label{gz}
\lim_{\delta\to 0+}\,\int_0^\delta \sup_{|\beta|\le r}\,|f_k'(\alpha+i\beta)|\,d\alpha=0,
\end{equation*}
uniformly in $k$. Let $g \in \Bes$ with $\lim_{|z|\to\infty} g(z)=0$.  
Then $f\in \mathcal{B}$, and 
\[
\lim_{k\to\infty}\,\|(f_k-f)g\|_{\mathcal{B}}=0.
\]
\end{lemma}

\subsection*{Shifts and density results for $\B$}

In order to establish the existence of a $\B$-calculus for many operators, we will need a dense subalgebra of $\Bq$ with certain properties.     To establish the density, we will make use of the shift operators $(T_\Bes(\tau))_{\tau\in\overline\C_+}$ on $\Bes$, which are defined  by the formula
\[
(T_\Bes(\tau)f)(z) := f(z+\tau), \qquad  f \in \Bes, \, \tau\in\overline\C_+, \, z\in\C_+.
\]
The following properties are established in \cite{BGT}: see Lemma 2.6 for statements (1)-(4),  the proof of Proposition 2.10 for (5) and (6), and Proposition 4.6 for (7).

\begin{prop}  \label{lem2.6}
\begin{enumerate}[\rm1.]
\item   For each $f \in \Bes$, 
\[
\|T_\Bes(\tau)f\|_\Bes \le \|f\|_\Bes, \quad \tau \in \overline\C_+.
\]
Moreover
\[
\lim_{\tau\in\overline\C_+, \tau\to0}\,\|T_\Bes(\tau)f-f\|_\Bes=0.
\]
\item  Let $-A_\Bes$ be the generator of the $C_0$-semigroup $(T_\Bes(t))_{t\ge0}$ on $\Bes$.  Then
\[
D(A_\Bes) = \{ f \in \Bes: f' \in \Bes\},  \quad A_\Bes f = -f'.
\]
\item The generator of the $C_0$-group $(T_\Bes(-is))_{s\in\R}$ is $iA_\Bes$.
\item $\sigma(A_\Bes) = \R_+$.
\item For each $f\in\Bq$,
\[
\lim_{t\to\infty} \|T_\Bes(t)f\|_\Bq  =  \lim_{t\to\infty} \Big\|f - A_\Bes \int_0^t T_\Bes(s)f \, ds \Big\|_{\Bq} = 0.
 \]
\item The closure of the range of $A_\Bes$ is $\Bq$.
\item The family $(T_\Bes(\tau))_{\tau\in\C_+}$ is a holomorphic $C_0$-semigroup of contractions.
\end{enumerate}
\end{prop}

As a preliminary exercise, we consider the rational functions in $\Bq$.   
Let $\Rat$ be the  linear span of $\{r_\l: \l\in\C_+\}$.   It follows from the resolvent identity and the continuity of the
map $\l \mapsto r_\l$ from $\C_+$ to $\LT$, that $\Rat$ is dense in the algebra of all rational functions which vanish at infinity and have no poles in $\C_+$, with respect to the HP-norm and hence with respect to the $\B$-norm and the $H^\infty$-norm.  

\begin{prop} \label{CCCA}  
The closure of $\Rat$ in $\B$, with respect to either the  $\B$-norm or the $H^\infty$-norm, is $\B \cap C_0(\overline\C_+)$. 
\end{prop}

The case of the $H^\infty$-norm in Proposition \ref{CCCA} is well-known (see \cite[Appendix F, Proposition F.3]{HaaseB}, for example).   The coincidence of the closures in the different norms is a special case of a more general result \cite[Theorem 4.4]{BGT2} for subspaces of $\B$ which are invariant under the horizontal shifts $(T_\B(\tau))_{\tau\ge0}$ considered above.   As stated in Proposition \ref{balg}, the Hille-Phillips algebra $\LT$ is neither closed nor dense in $\Bes$.   We do not know of any explicit description of the closure of $\LT$ in $\Bes$, but the closures in the $\B$-norm and in the $H^\infty$-norm coincide. 







We now identify a dense subalgebra $\G$ of $\Bq$ consisting of (restrictions of) entire functions of exponential type.  Various manipulations can be carried out initially for these functions, and then extended to all functions in $\Bq$ by density arguments.   Thus the subalgebra $\G$ plays a very important role in the construction and theory of the $\B$-calculus.

It is well-known that the map $f \mapsto f^b$ is an isometric isomorphism from $H^\infty(\C_+)$ onto
\[
H^\infty(\R) :=\left\{ g \in L^\infty(\R): \operatorname{supp}(\mathcal{F}^{-1}g) \subset \R_+ \right\}
\] 
where $\F$ is the (distributional) Fourier transform  (see \cite[Section II.1.5]{Havin}, for example).  Let $I$ be a closed subset of $\R_+$.  The spectral subspace $H^{\infty}_I$ of $H^\infty(\C_+)$ is defined to be
\[
H^{\infty}_I = \big\{f \in H^\infty(\C_+): \operatorname{supp} \big(\fti f^b\big) \subset I\big\}.
\]

\begin{prop} \label{gest}
Let $I := [\ep,\sigma]$ and $J$ be non-empty compact intervals in $(0,\infty)$, and $f \in H^\infty_I$, $g \in H^\infty_J$.  
\begin{enumerate}[\rm1.]
\item $f\in \Bes$ and 
\begin{equation} \label{Gest}
\|f\|_\Bes \le \|f\|_\infty \left( 1 + 4 \log\left(1 + \frac{\sigma}{\ep}\right) \right).
\end{equation}
\item $fg \in H^\infty_{I+J}$.
\end{enumerate}
\end{prop}

See \cite[Lemma 2.5]{BGT} and \cite[Lemma VI.4.7]{Katz} for proofs.   It follows that
\begin{equation} \label{Gdef}
\G := \bigcup \left\{ H^\infty_I: \text{$I \subset (0,\infty)$, $I$ compact} \right\}
\end{equation}
is a subalgebra of $\Bq$.   The following result is proved in \cite[Proposition 2.10]{BGT}.

\begin{thm} \label{Gdense}
The closure of the subalgebra $\G$ in $\Bes$ is $\Bq$.
\end{thm}

We outline the proof of this statement, for which we will use Arveson's spectral theory for representations of locally compact abelian groups on Banach spaces (see \cite{Arv}, for example) in the special case of $C_0$-groups of isometries. For a bounded $C_0$-group $(T(t))_{t\in\R}$ on a Banach space $X$, let $\operatorname{sp}_T(x)$ be the Arveson spectrum of an element $x \in X$ with respect to $T$.  For a closed subset $I$ of $\R$, let $X_T(I) =\{x \in X:  \operatorname{sp}_T(x) \subset I\}$.

We apply this theory to the $C_0$-group of shifts on $\operatorname{BUC}(\R)$:
\[
(S(t))(s) = (s-t), \qquad s,t \in \R, \; g \in \operatorname{BUC}(\R),
\]
and to the $C_0$-group of vertical shifts:
\[
G(t) := T_\Bes(-it), \qquad t\in\R,
\]
on $\Bq$, as in Proposition \ref{lem2.6}.   Let $K : \Bes \to \operatorname{BUC}(\R)$ be the isometric injection given by $Kf = f^b$.   Then
\[
S(t)K = KG(t), \qquad t\in\R.
\]
For $f\in\B$, it follows from this and the definition of the Arveson spectrum that $\operatorname{sp}_G(f) = \operatorname{sp}_S(f) 
=\operatorname{supp} (\F^{-1}f^b)$.   Hence $ H^\infty_I = \B_G(I)$ 
 for all compact subsets $I$ of $(0,\infty)$.  

The density of $\G$ in $\Bq$ follows from a combination of Proposition \ref{lem2.6} and a further abstract result given in Proposition \ref{dense} below.  It is a variant of the fact that the Arveson spectrum of a bounded $C_0$-group, with generator $A$, is $-i\sigma(A)$ \cite[Theorem 8.19]{Dav80}.  A direct proof is given in \cite[Proposition 2.8]{BGT}.

\begin{prop}\label{dense}
Let $(T(t))_{t \ge 0}$ be a bounded $C_0$-group on a Banach space $X$, with generator $A$, and assume that
the range of $A$ is dense in $X$. 
\begin{enumerate}[\rm1.]
\item The set $\bigcup \{X_T(I) : I \; {\rm compact}, 0\notin I \}$ is dense in $X$.
\item If, in addition, $\sigma(A)\subset i [0,\infty)$, then $\bigcup \{X_T(I) : I \; {\rm compact}, I \subset (0,\infty) \}$ is dense in $X$.
\end{enumerate}
\end{prop}

\subsection*{The space $\E$ and duality}

Let $\E$ be the space of all holomorphic functions $g$ on $\C_+$ such that
\[
\|g\|_{\E_0} := \sup_{\a>0} \,\, \a \int_\R |f'(\a+i\b)| \, d\b < \infty.
\]
If $g \in \E$, then $g(\infty) := \lim_{\Re z \to\infty} g(z)$ exists in $\C$.  A norm on $\E$ is defined by
\[
\|g\|_\E := |g(\infty)| + \|g\|_{\E_0}.
\]
Then $\E$ is a Banach space.

There is a (partial) duality between $\Bov$ and $\Bes$ given by
\begin{equation} \label{dualdef}
\langle g,f\rangle_\Bes:=\int_0^\infty \a\int_\R\,g'(\a-i\b) f'(\a+i\b)\,d\b \,d\a,  \qquad g\in \Bov, \quad f\in \Bes.
\end{equation}
It follows from the definitions of $\Bes$ and $\Bov$ that the integral exists and 
\begin{equation} \label{ebbd}
|\langle g,f \rangle_\Bes |\le \|g\|_{\Bov_0}  \|f\|_{\Bes_0}.
\end{equation}
For $z \in \C_+$, the function $r_z\in\E$ and the formula \eqref{Brepf} for $f \in\B$ can now be written as 
\[
f(z) = f(\infty) + \langle r_z, f \rangle_\B.
\]

The duality is only partial in the sense that $\Bes$ and $\Bov$ are not the dual or predual of each other with respect to this duality, and the constant functions in each space are annihilated in the duality.  The duality appeared in \cite[p.266]{V1}  (in slightly different form), and it was
 noted there that Green's formula on $\C_+$ transforms (\ref{dualdef}) into 
\begin{equation} \label{Green}
\langle g,f \rangle_\Bes = \frac{1}{4} \int_\R  g(-i\b)f(i\b) \, d\b
\end{equation}
for ``good functions".  For sample results about the interpretation of ``good functions", see \cite[Lemma 17]{Taib} and \cite[Proposition 2.2]{BGT3}.
Informally, suppose that $f(\infty)=0$ (without essential loss of generality), and that \eqref{Green} holds for $g = r_z$, where $z \in \C_+$.  Then \eqref{Brepf} becomes
\[
f(z) = \frac{1}{2\pi} \int_\R \frac{f(i\b)}{z-i\b} \, d\b,
\]
which is a Cauchy formula for $f(z)$.   However the integral here is not necessarily absolutely convergent, and the formula is not valid for some $f \in \Bes$.  

The duality plays a crucial role in the construction of the $\B$-calculus, via the following approximation lemma.   Theorem \ref{Gdense}, shows that, in order to prove the lemma, it suffices to consider the case when $f \in \G$, and then an inequality of Bohr's type is applicable.  See \cite[Lemma 3.19]{BGT} for details of the proof.

\begin{lemma}\label{Sconv}
Let $\mu \in M(\R_+)$ satisfy
\[
\quad \mu(\R_+)= 1, \quad  \int_{\R_+} t\,d|\mu|(t)  <\infty.
\]
Let $m = \lt\mu \in \LT$, and $f\in \Bes$.  For $\delta>0$, let
\[
 f_\delta(z)  := m(\delta z) f(z), \quad z \in \C_+.
\]
Then, for all $g \in \Bov$,
\[
\lim_{\delta\to0+} \langle g, f_\delta \rangle_\Bes = \langle g,f \rangle_\Bes.
\]
\end{lemma}


\subsection*{Definition of the $\Bes$-calculus}

We can now proceed towards a definition of a $\Bes$-calculus by replacing $z$ by an operator $A$ in Proposition \ref{BHP2}.   To make this effective for as many operators $A$ as possible, we will consider the integral in \eqref{Brepf} in the weak operator topology (see \eqref{fcdef} below).   For this we need the following assumption on the operator $A$.

Let $A$ be a closed densely defined operator on a Banach space $X$, and assume that the spectrum $\sigma(A)$ is contained in $\overline\C_+$ and 
\begin{equation} \label{8.1}
\sup_{\alpha >0}\, \alpha \int_{\mathbb R} |\langle (\alpha +i\beta + A)^{-2}x, x^* \rangle| \, d\beta <\infty
\end{equation}
for all $x \in X$ and $x^* \in X^*$.   In other words, the functions 
\begin{equation} \label{resbov}
g_{x,x^*} : z \mapsto \langle (z+A)^{-1}x, x^* \rangle
\end{equation}
belong to $\Bov$.   These assumptions were first considered by Gomilko \cite{Gom}, and  by Shi and Feng \cite{SF}, who independently proved that they imply that $-A$ is the generator of a bounded $C_0$-semigroup on $X$.  We shall call them the (GSF) {\it assumptions} on $A$.

By the Closed Graph Theorem, the supremum in \eqref{8.1} is bounded for $\|x\|\le1$ and $\|x^*\|\le1$.  Hence, if the (GSF) assumptions are satisfied there is a minimal constant $\gamma_A$ such that
\begin{equation} \label{gammaa}
\int_{\mathbb R} |\langle (\alpha +i\beta + A)^{-2}x, x^* \rangle| \, d\beta \le \frac{\gamma_A}{\a} \|x\|\,\|x^*\|, \qquad x \in X, \,x^*\in X^*.
\end{equation}

\begin{exas}
1.  Let $-A$ be the generator of a bounded $C_0$-semigroup $(T(t))_{t\ge0}$ on a Hilbert space $X$.
  Then $A$ satisfies the (GSF)  assumptions with $\gamma_A \le 2K_A^2$  \cite[Example 4.1(1)]{BGT}.

\noindent 2.  Let $A \in \Sect(\pi/2-)$.  Then $A$ satisfies the (GSF) assumptions, with $\gamma_A \le C M_A (\log M_A + 1)$, where $C$ is an absolute constant.   See Theorem \ref{sectB} below.

\noindent 3.  Let $A$ be the generator of the $C_0$-group of translations on $L^p(\R)$, where $1 \le p < \infty, \, p\ne2$.   Then $A$ does not satisfy the (GSF) assumptions \cite[Corollary 6.7]{BGT2}.
\end{exas}

Let $A$ satisfy the (GSF) assumptions, and let $f \in \Bes$.  Define
\begin{gather}  \label{fcdef}
\null\hskip-70pt\langle f(A)x, x^* \rangle :=  f(\infty) \langle x, x^* \rangle + \frac{2}{\pi} \langle g_{x,x^*}, f \rangle_\Bes\\
\null\hskip30pt =  f(\infty) \langle x, x^* \rangle - \frac{2}{\pi}  \int_0^\infty  \alpha \int_{\R}   \langle (\alpha -i\beta +A)^{-2}x, x^* \rangle  {f'(\alpha +i\beta)} \, d\beta\,d\alpha  \nonumber 
\end{gather}
for all $x \in X$ and $x^* \in X^*$.  It is easily seen that this defines a bounded linear mapping $f(A) : X \to X^{**}$, and that the linear mapping
\[
\Phi_A : \Bes \to L(X,X^{**}), \qquad f \mapsto f(A),
\]
is bounded.  Indeed, by \eqref{ebbd} and \eqref{gammaa},
\begin{equation*} \label{const}
\|f(A)\| \le |f(\infty)| + \gamma_A \|f\|_\Bq \le \gamma_A\|f\|_\Bes.
\end{equation*}

\begin{thm} \label{besc}
If $A$ satisfies the {\rm(GSF)} assumptions, then the map $\Phi_A : f \mapsto f(A)$ is a bounded algebra homomorphism from $\Bes$ into $\L(X)$, which extends the Hille-Phillips calculus.  Moreover, $\|\Phi_A\| \le \gamma_A$.   
\end{thm}

It is not immediately clear why \eqref{fcdef} defines an operator $f(A)$ which maps $X$ into $X$.  The proof of this has three steps which we sketch here.  

Firstly, consider $f \in \LT$.   Then the calculation at the beginning of this section is justifiable in the weak operator topology, and it shows that $f(A)x$ as in \eqref{fcdef} agrees with $f(A)x$ as defined in \eqref{HP}.   In particular, $f(A)$ maps $X$ into $X$.

Secondly, consider $f \in \G$.   For $\delta>0$, let
$$
e_\delta(z) := \frac{1-e^{-\delta z}}{\delta z}, \qquad z \in \C_+.
$$
Then $e_\delta \in \LT$,  $f e_\delta \in \LT$,  and Lemma \ref{Sconv} can be used to show that $\lim_{\delta\to0+} (f e_\delta)(A) = f(A)$ in the strong operator topology.  Thus $f(A)$ maps $X$ into $X$.

Finally, consider $f \in \Bes$.  By Theorem \ref{Gdense}, $f - f(\infty)$ is in the norm-closure of $\G$, so $f(A) - f(\infty)$ is the limit in operator norm of a sequence $(f_n(A))_{n\ge1}$ where each $f_n \in \G$.   Hence $f(A)$ maps $X$ into $X$.

For more detailed proofs, including the fact that $\Phi_A$ is an algebra homomorphism, see \cite[Lemma 4.2, Lemma 4.3 and Theorem 4.4]{BGT}.

There are converse results to Theorem \ref{besc} as follows.  
\begin{thm} \label{cseb}
Let $A$ be a closed densely defined operator on a Banach space $X$, and assume that $\sigma(A) \subseteq \overline\C_+$ and there exists a $\Bes$-calculus $\Phi$ for $A$.   Then $A$ satisfies the {\rm(GSF)} assumptions, and $\Phi=\Phi_A$ as defined in \eqref{fcdef}.
\end{thm}


The proof of necessity of the (GSF) assumptions uses the general method of constructing functions in $\B$ described in Proposition \ref{BHP2}, and the proof that $\Phi=\Phi_A$ can be obtained by following the steps in the proof of Theorem \ref{besc}.  Details are given in  \cite[Theorems 6.1 and 6.2]{BGT2}.

Theroems \ref{besc} and \ref{cseb} establish that a densely defined operator $A$ has a bounded functional calculus for $\Bes$ if and only if $A$ satisfies the (GSF) assumptions.   The unique functional calculus is $\Phi_A$ as above, and it extends the Hille-Phillips calculus.   We will say that $\Phi_A$ is the {\it $\Bes$-calculus} for $A$.

The $\Bes$-calculus is compatible with the sectorial calculus in various ways.  In particular if $A \in \Sect(\pi/2-)$, $f \in \Bes$ and $f(A)$ can be defined in the sectorial calculus (for example, if $A$ is injective), then the definitions agree.   There are also similar statements when $A$ satisfies the (GSF) conditions, $f$ is holomorphic on a sector larger than $\C_+$, and the restriction of $f$ to $\C_+$ belongs to $\Bes$.   The $\Bes$-calculus is compatible in similar ways with the half-plane calculus considered in \cite{BHM}.  We refer the reader to \cite[Section 4.3]{BGT} for details.

There is a Convergence Lemma for the $\Bes$-calculus, which can easily be deduced from Lemma \ref{lem9.1} (see \cite[Corollary 8.3]{BGT2}).  It was first proved by a different argument in \cite[Theorem 4.13]{BGT}.

\begin{thm} \label{cor8.3}
Let $A$ be an operator which satisfies the {\rm(GSF)} assumptions, and let $f_k$ and $f_0$ be as in Lemma \ref{lem9.1}.  Then $f_k(A) \to f_0(A)$ in the strong operator topology as $k\to\infty$.
\end{thm}

There is a Spectral Inclusion Theorem for the $\B$-calculus, as follows.   The proofs are simple variants of arguments used for the HP-calculus, involving Banach algebra techniques as in \cite[Section 16]{HP}.   The third statement uses the $F$-product as in \cite[Section A.III.4]{Nagel}.   See \cite[Theorem 4.17]{BGT} for the proofs of the other statements.

\begin{thm}  \label{thm9.3}
 Let $A$ be an operator which satisfies the {\rm(GSF)} assumptions.  Let $f \in \Bes$ and $\l \in \C$.
\begin{enumerate} [\rm1.]
\item  If $x \in D(A)$ and $Ax = \l x$, then $f(A)x = f(\l)x$.
\item  If $x^* \in D(A^*)$ and $A^*x^* = \l x^*$, then $f(A)^*x^* = f(\l)x^*$.
\item  If $(x_n)_{n\ge1}$ are unit vectors in $D(A)$ and $\lim_{n\to\infty}\|Ax_n - \l x_n\| = 0$, then $\lim_{n\to\infty}\|f(A)x_n - f(\l)x_n\| = 0$.
\item  If $\l \in \sigma(A)$ then $f(\l) \in \sigma(f(A))$.
\item If $A \in \Sect(\pi/2-)$, then $\sigma(f(A)) \cup \{f(\infty)\} = f(\sigma(A)) \cup \{f(\infty)\}$.
\end{enumerate}
\end{thm}

\subsection*{Applications}   
The $\B$-calculus has various applications, some of them involving fine estimates of norms of particular operators $f(A)$, and others providing new proofs and/or variants of known results in semigroup theory.   We give some examples of each type.  Recall that $\L\mu$ denotes the Laplace transform of a measure $\mu\in M(\R_+)$.

\noindent
1.  {\it Spectral Mapping Theorem}.   The Spectral Mapping Theorem does not hold for $C_0$-semigroups even when $X$ is a Hilbert space, so it does not hold for the $\B$-calculus.   However the following result shows that it holds for a subalgebra of $\B$. 

\begin{thm} \label{thm9.3b}
If $-A$ generates a bounded $C_0$-semigroup on a Hilbert space, and $f \in \mathcal B \cap C_0({\mathbb C}_+)$, then
\begin{equation} \label{smt}
\sigma(f(A)) \cup \{0\} = f (\sigma(A))\cup \{0\}.
\end{equation}
In particular,
if the Fourier transform of $\mu\in M(\R_+)$ belongs to $C_0(\R)$ then
\[
\sigma(\mathcal L \mu (A)) \cup \{0\}=\mathcal L \mu (\sigma(A))\cup \{0\}.
\]
\end{thm}

This result can be proved by using  Proposition \ref{CCCA} and a result of Stafney \cite[Lemma 6.1]{Stafney}.   See \cite[Theorem 7.3]{BGT2} for details.   In fact, the result holds for all operators $A$ satisfying the (GSF) assumptions.

\noindent
2.  {\it Bernstein functions and subordination}.  Let $g$ be a Bernstein function on $\C_+$.   One of the many equivalent definitions of Bernstein functions says that there is a vaguely continuous convolution semigroup of subprobabiility measures $(\mu_t)_{t\ge0}$ on $\R_+$ such that $\L\mu_t = e^{-tg}$ for all $t\ge0$.    See \cite{Schill} for further information about Bernstein functions.

Phillips \cite{Phillips} showed that, if $-A$ is the generator of a bounded $C_0$-semigroup on $X$, then the operators $S(t)$ defined by 
\[
S(t)x := \int_0^\infty e^{-sA} x \,d\mu_t(s), \qquad x \in X, \; t\ge0,
\]
form a bounded $C_0$-semigroup on $X$, which is said to be {\it subordinate} to $X$.   The negative generator of the subordinate semigroup is an operator  $g_{\text{Bern}}(A)$.    This operator  $g_{\text{Bern}}(A)$ coincides with the operator $g(A)$ defined in other functional calculi, such as the extended sectorial calculus if $A \in \Sect(\pi/2-)$ and $A$ is injective.

The following question was raised in \cite{KR}:

If $A \in \Sect(\theta)$ where $\theta \in [0,\pi/2)$, is $g_{\text{Bern}}(A) \in \Sect(\theta)$?

\noindent
In other words, if $-A$ generates a bounded holomorphic $C_0$-semigroup of angle $\theta$, is the subordinate semigroup generated by $-g_{\text{Bern}}(A)$ also holomorphic, of angle at least $\theta$?    

Following some partial answers, a positive answer to the full question was given in \cite{GT15},  followed by a second proof in \cite{BGTad}.  The $\Bes$-calculus can be used to give a third proof \cite[Proposition 3.8, Corollary 5.10]{BGT}.    A further proof, based on the $\D$- and $\H$-calculi, is briefly discussed at the end of Section \ref{DH-calc}.

\noindent
3.  {\it Short-time behaviour}.  It is a starting point of semigroup theory that one defines the generator $-A$ of a $C_0$-semigroup as the right-hand derivative (at $t=0$) of the semigroup $e^{-tA}$ on its natural domain.  The result below, proved in \cite[Theorem 8.9]{BGT}, shows that the negative exponential function can be replaced by many other functions.
The proof uses the $\Bes$-calculus and other properties, so it is not clear that the result holds for all operators satisfying the (GSF) assumptions.

\begin{thm}\label{Gen2}
Let $-A$ be the generator of either a bounded $C_0$-semigroup on a Hilbert space $X$, or a bounded holomorphic $C_0$-semigroup on a Banach space $X$.
Let $f\in D(A_{\mathcal B})$, so 
$f\in \mathcal B$ and $f'\in \mathcal{B}$.    Then, for all $x \in D(A)$, 
\begin{equation} \label{ftd}
\lim_{t\to 0+} t^{-1} (f(tA)x-f(0)x) =f'(0)Ax.
\end{equation}
Conversely, if $f'(0)\ne0$ and 
\[
\lim_{t\to 0+} t^{-1} (f(tA)x-f(0)x)
\]
exists for some $x \in X$, then $x \in D(A)$ and \eqref{ftd} holds.    Thus, if $f'(0)=-1$, then $-A$ coincides with the right-hand derivative of $f(\cdot A)$ at $0$, in the strong operator topology, on its natural domain.
\end{thm}

\noindent
4.  {\it Cayley transform problem on Hilbert spaces}.  Let $-A$ be the generator of a bounded $C_0$-semigroup on a Hilbert space, and $V(A) = (A-1)(A+\nolinebreak[4]1)^{-1}$.  If the semigroup is contractive, it is well-known that $V(A)$ is also a contraction.  The discussion in Section \ref{HP-calc} leaves open the specific question whether $V(A)$ must be power-bounded when the semigroup is bounded but not contractive..    It was shown in \cite{Go04} that $\|V(A)^n\|$ grows at most logarithmically in $n$.   The proof of this involved intricate estimates of Laguerre polynomials.

The $\Bes$-calculus shows that 
\[
\|V(A)^n\| \le C \|V^n\|_{\Bes}.
\]
It is shown in \cite[Lemma 3.7]{BGT} (with relatively simple estimates) that $\|V^n\|_\Bes = O(\log n)$ 
 and in \cite[Lemma 5.1]{BGT2} that $\|V^n\|_\Bes \asymp \log n$.
 It remains unknown whether the estimate $\|V(A)^n\| = O(\log n)$  is optimal for bounded semigroups on Hilbert spaces (see Section \ref{A-calc} for further discussion).   

\noindent
5.  {\it Inverse generator problem  on Hilbert spaces}.
This problem was discussed in Section \ref{HP-calc}.  In the case of Hilbert spaces, the answer is unknown.    Let $A$ have dense range on a Hilbert space, and assume  that $-A$ generates a bounded $C_0$-semigroup.   If the inverse generator problem on Hilbert spaces has a positive answer, then the $C_0$-semigroups $(e^{-tA^{-1}})_{t\ge0}$ must be bounded.   Thus a subsidiary question is whether all semigroups of this type on Hilbert spaces are bounded.   This question can be addressed by functional calculus.

Let $\varphi_t(z) = z(z+1)^{-1} e^{-t/z}$.   Then $\varphi_t \in \LT$.  The HP-norm and the $\Bes$-norm can be estimated as follows \cite[Section 5.3]{BGT2}.

\begin{prop} \label{phit}
Let $\varphi_t, \, t\ge0$ be as above.  Then
\[
\|\varphi_t\|_{\rm{HP}} \asymp t, \qquad
\|\varphi_t\|_\Bes \asymp \log t, \qquad t\to \infty.
\]
\end{prop}

Applying the $\Bes$-calculus and the estimate above for $\|\varphi_t\|_\B$ to estimate the norm of $\varphi_t(A)$  leads easily to the following result, a version of which was obtained in \cite[Theorem 2.2]{Zwart07} under stronger assumptions on $A$.   Although the logarithmic estimate in Proposition \ref{phit} is sharp, it does not follow that the estimate in Corollary \ref{CH2} is sharp.

\begin{cor} \label{CH2}
Let $-A$ be the generator of a bounded $C_0$-semigroup on a Hilbert space $X$, and assume that $A$ has a bounded inverse.   Then 
\[
\|e^{-tA^{-1}}\| = O(\log t), \qquad t\to\infty.
\] 
\end{cor}

\noindent
6.  {\it Exponentially stable semigroups}.   Let $-A$ be the generator of an exponentially stable $C_0$-semigroup on a Hilbert space, and $f \in H^\infty(\C_+)$.   In this situation, $f(A)$ can be defined as a closed operator via the half-plane calculus.  Under a weak assumption on the boundary function of $f$, the following result shows that $f(A) \in L(X)$. 

\begin{cor}\label{CH}
Let $-A$ be the generator of an exponentially stable $C_0$-semigroup on a Hilbert space $X$, so that
\[
\|e^{-tA}\|\le Me^{-\omega t},\qquad t\ge 0,
\]
for some $M,\omega>0$.  Let  $f\in H^\infty(\mathbb C_{+})$, let
\[
h(t) = \operatorname{ess\,sup}_{|s|\ge t} |f(is)|,
\]
and assume that 
\[
\int_0^\infty \frac{h(t)}{1+t} \, dt < \infty.
\]
Then $f(A) \in L(X)$ and 
\[
\|f(A)\|\le 6 M^2  \int_0^\infty \frac{h(t)}{\omega + t} \, dt.
\]
\end{cor}

The proof  uses the $\Bes$-calculus, applied to the function $g(z) := f(z+\omega)$, which belongs to $\Bq$.   For the details, see \cite[Lemma 3.3 and Corollary 5.6]{BGT}. The result had previously been proved by Schwenninger and Zwart \cite{SZ} in the special case when $|f(is)| \le (\log(|s|+e)^{-\a}$ for some $\a>1$.   It is valid on Banach spaces, provided that $A$ satisfies the (GSF) assumptions.

\section{$\D$- and $\H$-calculi for sectorial operators}  \label{DH-calc}

Let $A \in \Sect(\pi/2-)$, so that \eqref{R11} holds 
and $-A$ is the generator of a bounded holomorphic $C_0$-semigroup on $X$.  In \cite{V1}, Vitse established that $A$ has a bounded $\B$-calculus with an estimate  $\|f(A)\| \le 31 M_A^3 \|f\|_\B$.   An estimate of the form \eqref{4.17} was later obtained in \cite[Theorem 5.2]{Sch1}.  Both those papers used the sectorial functional calculus and Littlewood-Paley decompositions to obtain the estimates.  The estimate \eqref{4.16} was shown in \cite[Corollary 4.8]{BGT} (with $C = 8(2+\log 2)$), using methods of complex analysis and functional analysis.

\begin{thm} \label{sectB}
Let $A \in \Sect(\pi/2-)$.   There is an absolute constant $C$  such that
\begin{equation} \label{4.16}
\gamma_A \le C M_A (\log M_A+1).
\end{equation}
Hence $A$ has a $\Bes$-calculus, and, for all $f \in \Bes$,
\begin{equation}  \label{4.17}
\|f(A)\| \le C M_A (\log M_A+1) \|f\|_\B,
\end{equation}
and 
\[
f(A) = f(\infty) - \frac{2}{\pi} \int_0^\infty \a \int_\R (\a-i\b+A)^{-2} f'(\a+i\b) \,d\b \,d\a,
\]
where the integral converges in operator norm.
\end{thm}

For sectorial operators, it is possible to extend the $\Bes$-calculus to more functions in two ways. In the $\Bes$-calculus, the functions $f$ are holomorphic on the half-plane $\C_+$, but $\sigma(A)$ is contained in a smaller sector (including the origin).  Consequently it may be possible to define $f(A)$ as a bounded operator, for some functions $f$ which are holomorphic, but not bounded, on $\C_+$.    More significantly, one can define a bounded calculus for holomorphic functions which are defined only on subsectors of $\C_+$.   In this section we describe two related calculi, the $\D$-calculus for functions defined on $\C_+$, and the $\H$-calculus for functions on sectors.

\subsection*{The spaces  $\D_s$}   
For $s>-1$,  let $\mathcal{V}_s$ be the  Banach space of (equivalence classes of) measurable functions $g: \C_{+}\to \mathbb C$ such that the norm
\begin{align}\label{norm}
\|g\|_{\mathcal{V}_s}
&:= \int_0^\infty \alpha^s\int_\R \frac{|g(\alpha +i\beta)|}{(\alpha^2+\beta^2)^{(s+1)/2}}\,d\beta\,d\alpha  < \infty.
\end{align}
If $\sigma>s$, then
\begin{equation} \label{vst}
\mathcal{V}_s\subset \mathcal{V}_{\sigma} \qquad \text{and} \qquad  \|g\|_{\mathcal{V}_{\sigma}}\le \|g \|_{\mathcal{V}_s},
\qquad g\in \mathcal{V}_s.
\end{equation}

For $s>-1$, let $\D_s$ be the space of all holomorphic functions $f$ on $\C_+$ such that $f' \in \V_s$, with the seminorm
\[
\|f\|_{\D_{s,0}} := \|f'\|_{\V_s}.
\]
Functions in this class appeared in \cite[Corollary 4.2]{Al2014}, but an extended theory first appeared in \cite{BGT3}.

\begin{exas} \label{ex3.3}
 1.  Let $\l \in \C_+$, and consider $r_\l(z) := (\l+z)^{-1}$.  Then $r_\l \in \D_{s}$ for all $s>-1$, and $\|r_\l\|_{\D_{s,0}} \le C_s |\l|^{-1}$, where $C_s$ depends only on $s$.
	See \cite[Example 3.3]{BGT3}.
	
	\noindent
2.  Let $V(z) = (z-1)(z+1)^{-1}$.   For $s>0$, 
\[
\|V^n\|_{\D_{s,0}} \le 16 \big( B(s/2,1/2) + 2^{-s/2} \big), \qquad n\in\N,
\]
where $B$ is the beta function.  See \cite[Lemma 12.1]{BGT3}.   

\noindent
3. Let
\[
f_\nu(z):=z^{\nu} e^{-z},\quad z\in \C_{+},\quad \nu\ge0.
\]
Then $f_\nu\in \mathcal{D}_s$ if and only if $s>\nu$.
Moreover, if $s >\nu$, then
\begin{equation}\label{frac}
\|f_\nu\|_{\D_{s,0}}\le 2 B\left(\frac{s-\nu}{2},\frac{1}{2}\right)\Gamma(\nu+1).
\end{equation}
In particular, the function $e^{-z} \in \B$ but it is not in $\D_0$.   See \cite[Example 3.4]{BGT3}.

\noindent
4.   As noted in Example \ref{exB}(2), the function $\arccot$ is not in $\B$.  However it is in $\D_s$ for all $s>-1$.   The function $g(z):=\exp(\arccot z)$ is bounded on $\C_+$ and it is in $\D_s$ for all $s>-1$.  However the boundary function of $g$ is not continuous, so $g \notin \Bes$.   Moreover, $\lim_{\ep\to0+} g(\ep+i)$ does not exist.   See \cite[Examples 3.5 and 3.19]{BGT3}.
\end{exas}	

The following inclusions are easily seen from \eqref{vst} and  \cite[Proposition 3.15]{BGT3}.

\begin{prop} \label{BDs}
\begin{enumerate}[\rm1.]
\item Let $\sigma>s>-1$.   Then $\D_s \subset \D_\sigma$. 
\item Let $s > 0$.  Then $\Bes \subset \D_s$.
\end{enumerate}
\end{prop}


For $g \in \V_s$, let
\begin{equation} \label{qdef}
(Q_s g)(z):= - \frac{2^s}{\pi}\int_0^\infty \alpha^s\int_\R \frac{g(\alpha+i\beta)}{(\alpha-i\beta+z)^{s+1}} \, d\b\,d\alpha,
\quad z\in \C_{+} \cup \{0\}.
\end{equation}
Note that the definition of $Q_1$ agrees with the operator $Q$ defined in \eqref{defQ}, but the domain of $Q_1$ is larger than the domain of $Q$.  

The properties of $\V_s$, $\D_s$ and $Q_s$ stated in Proposition \ref{Furt} below generally correspond to the properties of $\W$, $\Bes$ and $Q$ given in Proposition \ref{BHP2}.   One variation is that $Q_s$ does not map all functions in $\V_s$ into $\D_s$ but Proposition \ref{Furt2} shows that $Q_s$ does map all functions in $\V_s$ into $\D_\sigma$ for all $\sigma>s$.   The representation formula \eqref{qnrep} was shown in \cite[Corollary 4.2]{Al2014}.   A proof of the first statement in Proposition \ref{Furt} is given in \cite[Proposition 3.7]{BGT3}, and the existence of the sectorial limits $f(0)$ and $f(\infty)$ and the representation formula \eqref{qnrep} follow as corollaries. 

A function $f$ on $\C_+$ has {\it sectorial limits} at $0$ and at $\infty$ if the following limits exist in $\C$ for all $\psi \in (0,\pi/2)$:
\[
f(0) := \lim_{|z|\to0,z\in\Sigma_\psi}   f(z), \qquad  f(\infty) := \lim_{|z|\to\infty,z\in\Sigma_\psi} f(z).
\] 
The function $Q_sg$ in \eqref{qdef} has sectorial limits $(Q_sg)(0)$ given by \eqref{qdef} for $z=0$, and $(Q_sg)(\infty)=0$.
 
\begin{prop}\label{Furt}  Let $s>-1$.
\begin{enumerate}[\rm1.]
\item Let $g$ be a holomorphic function in $\V_s$.  Then $Q_sg \in \D_s$, and $(Q_sg)'=g$.
\item Let $f\in\D_s$.  Then $f'\in\V_s$, $f$ has sectorial limits at $0$ and $\infty$, and 
\begin{align} \label{qnrep}
f(z)&=f(\infty)+(Q_s f')(z)\\
&=f(\infty)- \frac{2^s}{\pi}\int_0^\infty \alpha^s\int_\R
\frac{f'(\alpha+i\beta)}{(\alpha-i\beta+z)^{s+1}} \, d\b\, d\alpha.\notag
\end{align}
\item There exists $g_s\in\V_s$ such that $Q_sg_s \notin \D_s$.
\end{enumerate}
\end{prop}

There is a norm on $\D_s$ given by
\[
\|f\|_{\D_s} := |f(\infty)| + \|f\|_{\D_{s,0}}.
\]
With this norm, $\mathcal{D}_s$ is a Banach space \cite[Corollary 3.11]{BGT3}, and the inclusions in Proposition \ref{BDs} are continuous.   On the other hand, $\D_s$ is not an algebra,  but $\D_\infty:= \bigcup_{s>-1} \D_s$ is an algebra \cite[Lemma 3.21]{BGT3}. 
  
Some functions  in $\D_s$ are unbounded; $\arccot$ is an example.   Let $\D_s^\infty := \mathcal{D}_s \cap H^\infty(\C_+)$.  There is a norm on $\D_s^\infty$ defined by
\[
\|f\|_{\D_s^\infty} = \|f\|_\infty + \|f\|_{\D_{s,0}}.
\]
With this norm, $\D_s^\infty$ is a Banach algebra.

The following results are proved in \cite[Proposition 3.6]{BGT3}.

\begin{prop}\label{Furt2}
Let $\sigma > s > -1$.   
\begin{enumerate}[\rm1.]
\item The restriction of $Q_\sigma$ to $\mathcal{V}_s$ is in $L(\mathcal{V}_s,\D_s)$.
\item $Q_s \in L (\mathcal{V}_s, \mathcal{D}_{\sigma})$.
\end{enumerate} 
\end{prop}
	
The representation formula \eqref{qnrep} is used to prove the second and third statements in the following proposition (see \cite[Lemmas 3.17 and 3.22, Corollary 5.5] {BGT3}).

\begin{prop} \label{lem3.17}
Let $f \in \D_s$, where $s>-1$.  
\begin{enumerate}[\rm1.]
\item The functions $\tilde f(z):= f(1/z)$ and $f_t(z):=f(tz),\, t>0$, are in $\D_s$, and $\|\tilde f\|_{\D_{s,0}} = \|f\|_{\D_{s,0}}$, $\|f_t\|_{\D_s} = \|f\|_{\D_s}$. 
\item For $\tau \in \C_+$, let $T(\tau)f(z) = f(z+\tau)$.   Then $T(\tau)f \in \D_s$. Moreover $(T(\tau))_{\tau\in\C_+}$ form a bounded holomorphic $C_0$-semigroup on $\D_s$.
\item For $n\in\N$, the function $z^nf^{(n)}(z) \in \D_{s+n}$.
\end{enumerate}
\end{prop}

Note that $\D_s$ is not invariant under vertical shifts (see \cite[Example 3.19]{BGT3} or Example \ref{ex3.3}(4) above).

\subsection*{The spaces  $\H_\psi$}

Here we introduce spaces $\H_\psi$ of holomorphic functions on sectors $\Sigma_\psi$.   Later in this section we will construct a bounded functional calculus for these spaces, for operators $A \in \Sect(\theta)$ where $0 \le \theta < \psi < \pi$.

Let $\psi \in (0,\pi)$.  The Hardy space $H^1(\Sigma_\psi)$ on the sector $\Sigma_\psi$ is defined to be the space of all holomorphic  functions $g $ on $\Sigma_\psi$ such that
\begin{equation}\label{CC2}
\|g\|_{H^1(\Sigma_\psi)} := \sup_{|\varphi|<\psi}
\int_0^\infty \bigl(|g(te^{i\varphi})|+
|g(te^{-i\varphi})| \bigr) \,dt <\infty.
\end{equation}
For $\psi = \pi/2$, this definition may appear to differ from the classical definition of a Hardy space on a half-plane, but actually the two definitions coincide (a proof is given in \cite[Corollary 4.3]{BGT3}).

We consider the space  $\H_\psi$ of all holomorphic functions $f$ on $\Sigma_\psi$ such that $ f' \in H^1(\Sigma_\psi)$.
These spaces are sometimes called {\it Hardy-Sobolev spaces}.   

It was shown in \cite[Proposition 2.4]{BGT} that $\H_{\pi/2} \subset \B$ and $\lim_{|z|\to\infty, z\in\C_+} f(z)$ exists in $\C$ for all $f \in \H_{\pi/2}$ (see also Theorem \ref{lem4.13}(1) below).   The first two statements in Proposition \ref{thm4.8} follow immediately for $\psi=\pi/2$.   Since the map $f(z) \mapsto f(z^{2\psi/\pi})$ is a linear isomorphism of $\H_\psi$ onto $\H_{\pi/2}$, those statements also hold for all $\psi\in(0,\pi)$.   The third statement is easily seen.

\begin{prop} \label{thm4.8}
Let $f \in \H_\psi,\, \psi\in (0,\pi)$.
\begin{enumerate}[\rm1.]
\item The function $f$ extends to a continuous bounded function on $\overline{\Sigma}_{\psi}$.
\item The limit
\[
f(\infty) := \lim_{|z|\to\infty,  z \in\Sigma_\psi} f(z)
\]
exists in $\C$.
\item For all $z \in\C_+$, 
\begin{equation*}\label{bound_for_sup}
|f(z)| \le |f(\infty)| + \|f'\|_{H^1(\Sigma_\psi)}.
\end{equation*}
\end{enumerate}
\end{prop}

Consequently, there is a norm on $\H_\psi$ defined by
\[
\|f\|_{\H_\psi} := \|f\|_{H^\infty(\Sigma_\psi)} + \|f'\|_{H^1(\Sigma_\psi)}.
\]

\begin{prop} \label{lem4.9}
For any $\psi \in (0,\pi)$, the space $\H_\psi$ is a Banach algebra.   Moreover, for $\psi_1,\psi_2 \in (0,\pi)$, the map $f(z) \mapsto f(z^{\psi_1/\psi_2})$ is an isometric algebra homomorphism from $\H_{\psi_1}$ onto $\H_{\psi_2}$.
\end{prop}

\begin{exas} \label{ex4.10}
1.  Let $\psi \in (0,\pi)$, $\varphi \in (0, \pi-\psi)$ and $\l \in \Sigma_{\varphi}$.   Then $r_\l \in \H_\psi$, and 
\[
\|r_\l\|_{\H_\psi} \le \frac{C_{\varphi,\psi}}{|\l|},
\]
where $C_{\varphi,\psi}$ depends only on $\varphi$ and $\psi$.   See \cite[Example 4.10(1)]{BGT3}.

\noindent
2.   Let  $\gamma \in (0,\pi/(2\psi))$ and $e_{\gamma,\l}(z) := e^{-\l z^\gamma}, \,z \in \Sigma_\psi$.
Then $e_{\gamma,\l} \in \H_\psi$ and 
\begin{multline*}
\|e'_{\gamma,\l}\|_{H^1(\Sigma_\psi)} = \|e'_{1,\l} \|_{H^1(\Sigma_{\gamma\psi})}\le \int_{\partial\Sigma_{\gamma\psi}} |\l| e^{-\Re\l z} \, |dz|  \\
\le \frac{1}{\cos(\varphi+\gamma\psi)} + \frac{1}{\cos (\varphi-\gamma\psi)}.
\end{multline*}
\end{exas}

The estimates for resolvents of Bernstein functions in the next proposition illustrate the value of the spaces $\D_s$ and $\H_\psi$, and they are applied at the end of this section.   They are proved in  \cite[Lemma 3.20 and Corollary 4.16]{BGT3}.

\begin{prop}  \label{cor4.16}
Let $g$ be a Bernstein function. 
\begin{enumerate}[\rm1.]
\item  Let  $\l \in \C_+$ and $f(z) = (\l + g(z))^{-1}$ for $z \in \C_+$. Then $f \in \D_s$ for $s>2$ and
\[
\|f\|_{\D_s} \le \frac{C_s}{|\l|},
\]
where $C_s$ is a constant depending only on $s$.
\item Let $\psi\in (0,\pi/2)$, $\varphi \in (\psi,\pi)$ and $\lambda \in \Sigma_{\pi-\varphi}$. Let $f(z) = (\l + g(z))^{-1}$ for $z \in \Sigma_\psi$.   Then $f \in \H_\psi$ and
\begin{equation}\label{bernst}
\|f\|_{\mathcal \H_\psi} \le
2\left(\frac{1}{\sin(\min(\varphi,\pi/2))}+\frac{2}{\cos\psi \sin^2((\varphi-\psi)/2)}\right)\frac{1}{|\lambda|}.
\end{equation}
\end{enumerate}
\end{prop}

We now lay out relations between $\B$, $\D_s$ and $\H_\psi$ which are proved in \cite[Theorem 4.12 and Lemma 4.13]{BGT3}.

\begin{thm} \label{lem4.13}
\begin{enumerate}[\rm1.]
\item  If $f \in \H_{\pi/2}$, then there exists $g \in L^1(\R_+)$ such that $f = f(\infty) + \mathcal{L}g$.  Hence $\H_{\pi/2} \subset \LT \subset \Bes$, and the inclusions are continuous.  
\item If $f \in \H_{\pi/2}$, then $f \in \D_s$ for all $s>-1$, and the inclusions are continuous.
\item  If $f \in \D_s$ for some $s>-1$, then (the restriction of) $f$ belongs to $\H_\psi$ for every $\psi \in (0,\pi/2)$, and the restriction maps are continuous.
\end{enumerate}
\end{thm}

The representation formula \eqref{qnrep} for $\D_s$ in Proposition \ref{Furt}(2) can be transformed into the formula \eqref{4.34} below for functions in $\H_\psi$, by means of Proposition \ref{lem4.9}, Theorem \ref{lem4.13} and Proposition \ref{thm4.8}.

\begin{prop} \label{prop4.17}
Let $f \in \H_\psi$, where $\psi \in (0,\pi)$.  Let $\nu := 2\psi/\pi$ and
\[
f_\nu(z) := f(z^\nu), \qquad z \in \C_+.
\]
Then
\begin{equation} \label{4.34}
f(z) = f(\infty) - \frac{1}{\pi} \int_0^\infty \int_\R \frac{f_\nu'(\a+i\b)}{\a-i\b+z^{1/\nu}} \,d\b\,d\a, \qquad z \in \Sigma_\psi\cup\{0\}.
\end{equation}
\end{prop}

The formula \eqref{4.34} can be recast by replacing the double (area) integral with a  line integral involving boundary values of $f'$ instead of scalings of $f'$, and offering a kernel of a different type which may sometimes be more useful.
The reproducing formula \eqref{arctan} for functions in $\H_\psi$ was formulated and proved in \cite[Proposition  4.19]{BGT3}, but it resembles, and was inspired by, \cite[Lemma 7.4]{Boy}.  The $\arccot$ function is defined in \eqref{arccot}.

\begin{prop}\label{L7.4}
Let $f \in \mathcal \H_\psi$, $\psi\in(0,\pi)$.   Let $\gamma:=\pi/(2\psi)$, and
\begin{equation}\label{fpsi}
 f_\psi(t):=\frac{f(e^{i\psi}t)+f(e^{-i\psi}t)}{2}, \quad t>0.
\end{equation}
Then
\begin{equation}\label{arctan}
f(z)=f(\infty)-\frac{2}{\pi}\int_0^\infty
f_\psi'(t)\, {\arccot}(z^\gamma/t^{\gamma})\,dt,\qquad z\in \Sigma_\psi\cup \{0\}.
\end{equation}
\end{prop}

\subsection*{Density of rational functions}
Let $\RR(\C_+)$ be the linear span of $\{r_\l : \l \in \C_+\}$ and the constant functions.   The first statement in the following theorem plays an important role in the $\D$- and $\H$-calculi.

\begin{thm} \label{thm5.1}
\begin{enumerate}[\rm1.] 
\item The space $\RR(\C_+)$ is dense in $\D_s$ for $s>-1$.  
\item The space $\Bes$ is not dense in $\D_s^\infty$ for $s>0$.
\end{enumerate}
\end{thm}

We outline the proof of the first statement in Theorem \ref{thm5.1} in the case when $s \in (-1,0)$.   For $f \in \D_s$, let
\[
R_f(\l) := -\frac{1}{\pi} f'(\l) r_{\overline\l}, \qquad \l\in\C_+.
\]
Then,  $R_f$ is a continuous function from $\C_+$ to $\D_s$, and by Example \ref{ex3.3}(1) it is Bochner integrable with respect to area measure $S$ on $\C_+$.   Since the point evaluations are continuous on $\D_s$, it follows that
\[
Q_0 f' = \int_{\C_+} R_f(\l) \,dS(\l).
\]
This is the Bochner integral of a continuous function, so $Q_0f'$ belongs to the closure in $\D_s$ of the span of the integrand which is contained in the closure of $\Rat$.   Since $f = f(\infty) + Q_0f'$, it follows that $f$ is in the closure of $\RR(\C_+)$ in $\D_s$. 

A similar argument may be used to prove the first statement  for $s\ge0$.  

The second statement in Theorem \ref{thm5.1} follows from the fact that any function in the closure of $\Bes$ has a continuous extension to $i\R$. See \cite[Theorem 5.1, Corollary 5.2]{BGT3} for details.

There is a similar result for the spaces $\H_\psi$, where $\psi \in (0,\pi)$.  Let $\mathcal{R}_0(\Sigma_\psi)$ be the linear span of $\{r_\l : \l\in\Sigma_{\pi-\psi}\}$, and $\mathcal{R}_\H(\Sigma_\psi)$ be the closure of $\mathcal{R}_0(\Sigma_\psi)$ in $\H_\psi$.  See \cite[Theorem 5.10]{BGT3} for a proof of the following density result.

\begin{thm} \label{thm5.10}
For $\psi\in (0,\pi)$,
\[
\mathcal{R}_\H(\Sigma_\psi) = \left\{f \in \H_\psi : f(\infty) = 0\right\}.
\]
\end{thm}

\subsection*{Convergence Lemmas}  To obtain Convergence Lemmas for functions in $\D_s$ and $\H_\psi$, we need to consider functions of the form  $f(z^\gamma)$, where $\gamma\in(0,1)$.   See \cite[Corollary 4.14 and Lemmas 6.1, 6.2]{BGT3} for proofs.

\begin{lemma} \label{cor4.14}
Let $s>-1$, $f \in \D_s$ and $\gamma\in(0,1)$.  Let $f_\gamma(z) = f(z^\gamma), \, z\in\C_+$.  Then $f_\gamma \in \D_\sigma^\infty \cap \H_{\pi/2}$ for all $\sigma>-1$.   Moreover, the maps $f \mapsto f_\gamma$ are continuous from $\D_s$ into $\D_\sigma^\infty$ and from $\D_s$ into $\H_\psi$.
\end{lemma}

Now we state simultaneously Convergence Lemmas for $\D_s$ and $\H_\psi$.

\begin{lemma}[Convergence Lemmas]\label{lem6.1}
Let $\mathfrak A$ denote either $\D_s$ for some $s>-1$ or $\H_\psi$ for some $\psi \in (0,\pi)$.   
Let  $(f_k)_{k=1}^\infty\subset \frak A$ be such that
\begin{equation*}\label{Bound50A}
\sup_{k \ge 1}\|f_k\|_{\frak A}<\infty,
\end{equation*}
and for every $z\in \C_{+}$ there exists
\[
f(z):=\lim_{k\to\infty}\,f_k(z).
\]
Let $g\in \frak A$  satisfy
\[
g(0)=g(\infty)=0.
\]
Let $\gamma\in (0,1)$ and
\[
f_{k,\gamma}(z) := f_k(z^\gamma), \quad g_\gamma(z) := g(z^\gamma), \qquad z \in \C_+.
\]
Then
\[
\lim_{k\to\infty}\,\|(f_{k,\gamma}-f_{\gamma})g_\gamma\|_{\frak A}=0.
\]
\end{lemma}

\subsection*{The $\D$- and $\H$-calculi}

Let $A \in \Sect(\pi/2-)$.  To define the $\D$-calculus, we adapt the formula in \eqref{qnrep}, replacing $z$ by $A$.  If $f \in \D_s$ for some $s>-1$, define
\[
f_{\mathcal{D}}(A) := f(\infty) - \frac{2^s}{\pi}\int_0^\infty \a^s \int_\R f'(\a+i\b) (\a-i\b+A)^{-(s+1)} \,d\b\,d\a.
\]
This integral is absolutely convergent in operator norm.  If $f \in \Bes$ and $s=1$, this definition agrees with the formula in \eqref{fcdef}.

\begin{thm}  \label{thm7.4}
Let $A \in \Sect(\pi/2-)$.  
\begin{enumerate}[\rm1.]
\item If $A$ is injective, then $f_\D(A)$ as defined above coincides with $f(A)$ as defined in the sectorial functional calculus.
\item If $f \in \D_\infty$, then the definition of  $f_{\mathcal{D}}(A)$  does not depend on  $s$.
\item If $f \in \B \cap \D_\infty$, then $f_\D(A) = \Phi_A(f)$.
\item Let $\l \in \C_+$.  Then $(r_\l)_{\D}(A) = (\l+A)^{-1}$.
\item The map $\Psi_A : f \mapsto f_\D(A)$ is an algebra homomorphism from $\D_\infty$ to $L(X)$.
\item For each $s>-1$, the map $f \mapsto f_\D(A)$ is bounded from $\D_s$ to $L(X)$.
\end{enumerate}
Moreover, $\Psi_A$ is the unique linear map from $\D_\infty$ to $L(X)$ satisfying (4) and (6).
\end{thm}

Proofs of the statements in Theorem \ref{thm7.4} may be found in \cite[Theorems 7.4 and 7.6]{BGT3}.   The map $\Psi_A$ is called the {\it $\D$-calculus} for $A$.   An estimate for the boundedness of $\Psi_A$ on $\D_s$ is:
\begin{equation} \label{eq7.7}
\|f(A)\| \le |f(\infty)| + \frac{2M_A^{\lceil s \rceil +1}}{\pi} \|f\|_{\D_{s,0}}.
\end{equation}

Now suppose that $A\in\Sect(\theta)$, where $\theta \in [0,\pi)$, $\psi \in (\theta,\pi)$,  $\gamma:= \pi/(2\psi)$ and $f_{1/\gamma}(z) := f(z^{1/\gamma})$.  Then $f_{1/\gamma} \in \H_{\pi/2} \subset \D_0$, by Proposition \ref{lem4.9} and  Theorem \ref{lem4.13}(2).   Let
\begin{equation} \label{defH}
f_\H(A):= f(\infty) - \frac{1}{\pi} \int_0^\infty \int_\R f_{1/\gamma}'(\a+i\b) (\a-i\b+A^\gamma)^{-1} \,d\b\,d\a,
\end{equation}
where $A^\gamma$ is the fractional power of $A$.
This definition does not depend on the particular choice of $\psi$ and the consequent value of $\gamma$.   Moreover it agrees with $f(A)$ as defined by the sectorial functional calculus (assuming that $A$ is injective), and it agrees with $f_\D(A)$  (assuming that $\theta<\pi/2$ and $f \in \D_\infty$).   See \cite[Proposition 8.1]{BGT3}.   

Now we define the $\H$-calculus for $A$, as in \cite[Theorem 8.2]{BGT3}.

\begin{thm}  \label{thm8.2}
Let $A \in \Sect(\theta)$, where $\theta \in (0,\pi)$.   For any $\psi \in  (\theta, \pi)$ the formula \eqref{defH}
defines a bounded algebra homomorphism:
\begin{eqnarray*}
\Upsilon_A : \mathcal \H_\psi \mapsto L(X), \qquad \Upsilon_A(f)=f_\H(A).
\end{eqnarray*}
The homomorphism $\Upsilon_A$ satisfies
 $\Upsilon_A(r_\l)=(\l+A)^{-1}$ for all $\l \in \Sigma_{\pi-\psi}$, and it is the unique homomorphism with these properties.
\end{thm}

\begin{remark}
If $A$ is any operator which admits a $\D$-calculus satisfying the properties (4) and (6) of Theorem \ref{thm7.4}, then it follows from Example \ref{ex3.3}(1) that $A\in\Sect(\pi/2-)$.   Similarly, it follows from Example \ref{ex4.10} that if $A$ admits an $\H_\psi$-calculus as in Theorem \ref{thm8.2}, then $A \in \Sect(\theta)$ for some $\theta \in [0,\psi)$.
\end{remark}

Theorem \ref{Sovp} provides an alternative to the formula \eqref{defH}, based on \eqref{arctan} and an estimate 
\[
\|\arccot(A^\gamma)\|  \le \frac{\pi}{2} M_\psi(A)\|f\|_{\H_\psi}.
\]
There is a proof of the theorem in \cite[Lemma 8.4, Theorem 8.6]{BGT3}, followed by a discussion of relations between the theorem  and \cite{Boy}.
 
\begin{thm}\label{Sovp}
Let $A \in \Sect(\theta)$ and $\gamma=\pi/(2\psi)$, where $0 \le \theta < \psi < \pi$.  If $f \in \H_\psi$, and $f_\psi$ is given by \eqref{fpsi}, then
\begin{equation*}\label{formula_a}
f_{\mathcal H}(A)= f(\infty)-\frac{2}{\pi}\int_0^\infty
f_\psi'(t) \,\arccot(A^\gamma/t^{\gamma})\,dt,
\end{equation*}
where the integral converges in the operator norm topology, and
\begin{equation*}\label{NewEs1}
\|f_{\mathcal H}(A)\|
 \le |f(\infty)| + \frac{M_\psi(A)}{2} \|f'\|_{H^1(\Sigma_\psi)} \le {M_\psi(A)} \|f\|_{\H_\psi}.
\end{equation*} 
\end{thm}

\begin{rem}\label{rem_h}
If $\psi<\pi/2$ and $M_\psi(A)=1$ in \eqref{NewEs1}, that is, if $-A$ generates a holomorphic $C_0$-semigroup which is contractive on $\Sigma_{(\pi/2)-\psi}$,  then the $\mathcal H_\psi$-calculus is contractive.   We are not aware of such estimates in constructions of other calculi in the literature.
\end{rem}

\subsection*{Convergence Lemmas and Spectral Mapping Theorems}
There are Convergence Lemmas and Spectral Mapping Theorems for both the $\D$-calculus and the $\H$-calculus.
In the statements below, $f(A)$ may be interpreted as being either $f_\D(A)$ or $f_\H(A)$, as appropriate.  The Convergence Lemma in Theorem \ref{thm9.1} follows easily from Lemma \ref{lem6.1} applied with a suitable value of $\gamma$ and $g(z) = z(1+z)^{-2}$, for example.   See \cite[Section 6.1]{BGT3} for details.

\begin{thm}\label{thm9.1}
Let $A \in \Sect(\pi/2-)$, with dense range.  Let $f_k$ and $f$ be as in Lemma \ref{lem6.1}.   Then $f_k(A) \to f(A)$ in the strong operator topology as $k\to\infty$.
\end{thm} 


The Spectral Mapping Theorem is proved in a similar way to Theorem \ref{thm9.3}.

\begin{thm}  \label{thm9.3c}
 Let $A \in \Sect(\theta)$ where either (a) $\theta < \pi/2$ and $f \in\mathcal D_\infty$, or (b)  $\theta < \psi < \pi$ and $f \in \H_\psi$.  Let $\l \in \C$.
\begin{enumerate} [\rm1.]
\item  If $x \in D(A)$ and $Ax = \l x$, then $f(A)x = f(\l)x$.
\item  If $x^* \in D(A^*)$ and $A^*x^* = \l x^*$, then $f(A)^*x^* = f(\l)x^*$.
\item  If $(x_n)_{n\ge1}$ are unit vectors in $D(A)$ and $\lim_{n\to\infty}\|Ax_n - \l x_n\| = 0$, then $\lim_{n\to\infty}\|f(A)x_n - f(\l)x_n\| = 0$.
\item  One has $\sigma(f(A)) \cup \{f(\infty)\} = f(\sigma(A)) \cup \{f(\infty)\}$.
\end{enumerate}
\end{thm}

\subsection*{Applications}

The $\D$-calculus provides explicit norm-estimates for various known facts about an operator $A \in\Sect(\pi/2-)$, as the bounded holomorphic semigroup generated by $-A$ is $(e^{-tA})_{t\ge0}$ where $e^{-tA}$ is defined by any of the functional calculi.  The following corollary exhibits some estimates.  The power-boundedness of $V(A)$  was first proved in \cite{Cr93} and \cite{Palencia} using different methods.  The general form of the second and third statements are now classical.   The results follow from estimates in Examples \ref{ex3.3} and \eqref{eq7.7}, together with some compatibility arguments and Proposition \ref{lem3.17}(1).  See \cite[Corollary 10.1]{BGT3} for details.

\begin{cor} \label{cor10.1}
Let $A \in \Sect(\pi/2 -)$.
\begin{enumerate}[\rm1.]
\item
For all $n\in\N$, 
\[
\|V(A)^n\|\le 1+32(1+(\sqrt2\pi)^{-1})M_A^2.
\]
\item For all $\nu\ge0$ and $t>0$,
\[
\|A^\nu e^{-tA}\| \le 2^{\nu+2} t^{-\nu} \Gamma(\nu+1) M_A^{\lceil\nu\rceil+2}. 
\]
\item
If $A$ has dense range, then $A^{-1}$ generates
a bounded holomorphic $C_0$-semigroup $(e^{-tA^{-1}})_{t \ge 0}$ satisfying
\[
\|A^{-\nu} e^{-tA^{-1}}\| \le 2^{\nu+2} t^{-\nu} \Gamma(\nu+1) M_A^{\lceil\nu\rceil+2}
\]
for all $\nu\ge0$ and $t>0$.
\end{enumerate}
\end{cor}

The $\H$-calculus also provides norm-estimates, such as the following which can be deduced from Example \ref{ex4.10}(2)  and \eqref{NewEs1}.

\begin{cor} \label{better}
Let $A \in \Sect(\theta)$,  $\theta\in (0,\pi)$, and $\gamma  \in (0,\pi/(2\theta))$.  Then $(e^{-tA^\gamma})_{t\ge0}$ is a bounded holomorphic $C_0$-semigroup of angle $(\pi/2) - \gamma\theta$.  More precisely, if $\psi \in (\theta,\pi/(2\gamma))$ and $\lambda = |\l| e^{i\varphi} \in \Sigma_{(\pi/2)-\gamma\psi}$, then
\begin{equation} \label{est_gamma}
\|e^{-\lambda A^{\gamma}}\| \le \frac{1}{2}\left(\frac{1}{\cos(\gamma\psi+\varphi)}+\frac{1}{\cos(\gamma\psi-\varphi)}\right) M_{\gamma\psi}(A).
\end{equation}
\end{cor}

Another application of the $\D$- and $\H$-calculi provides a new proof that the question raised in Application 2 of Section \ref{B-calc} has a positive answer.  Let $g$ be a Bernstein function, and $A \in \Sect(\theta)$ where $\theta \in [0,\pi/2)$.  Applying the $\D$-calculus and Proposition \ref{cor4.16}(1) to functions of the form $(\l+g)^{-1}$ establishes that the operator $g_{\text{Bern}}(A) \in \Sect(\pi/2-)$.   Applying the $\H$-calculus and the estimate in \eqref{bernst} shows that $g_{\text{Bern}}(A) \in  \Sect(\theta)$.  See \cite[Theorem 10.5]{BGT3}.

\section{Further extensions on Hilbert spaces}  \label{A-calc}

In \cite{Pel}, Peller indicated two extensions of the Besov calculus for power-bounded operators on Hilbert spaces to larger classes of functions.   Here we briefly describe corresponding extensions in the case of bounded $C_0$-semigroups on Hilbert spaces.  
To put this in context, we need to explain the relationship between the space $\B$ of so-called analytic Besov functions and a conventional Besov space.   Briefly, the reason is that the space $\Bq$ is isometrically isomorphic to a space $\B_{\text{dyad}}$ which is a homogeneous Besov space of functions on $\R$ defined by means of Littlewood-Paley decompositions.    
This is described in the first subsection, and it is followed by outlines of results from \cite{ALM} and \cite{White}, which provide sharper extimates for the norm of $f(A)$ than the $\B$-norm does.  Although \cite{White} predated \cite{ALM} by three decades, it was accessible to other researchers only very recently, and we discuss \cite{ALM} before \cite{White} for presentational reasons.     Finally we discuss the effect of these results on the Cayley transform question for Hilbert spaces.  

\subsection*{$\B$ as a Besov space}

It is a standard fact that $H^\infty(\C_+)$ is isometrically isomorphic to  the space 
\[
H^\infty(\R) := \left\{g \in L^\infty(\R): \operatorname{supp} \mathcal{F}^{-1}g \subset \R_+ \right\},
\]
with the $\sup$-norm, where $\mathcal{F}^{-1}$ denotes the inverse (distributional) Fourier transform.   The isomorphism is obtained by passage to the boundary function, and its inverse is the standard extension via Poisson integrals::
\begin{align*} 
f \in H^\infty(\C_+) &\mapsto f^b \in H^\infty(\R),\\
g\in H^\infty(\R) &\mapsto \mathcal{P}g \in H^\infty(\C_+) \\
(\mathcal{P}g)(t+is) &:= \frac{1}{\pi} \int_\R \frac{tg(y)}{t^2 + (y-s)^2} \, dy, \quad t>0, s\in\R.
\end{align*}
The isomorphism $f \mapsto f^b$ maps $\Bq$ onto the space $\B_{\text{dyad}} := H^\infty(\R) \cap B^0_{\infty,1}(\R)$, where $B^0_{\infty,1}(\R)$ is a conventional homogeneous Besov space, defined in terms of Littlewood-Paley decomposition with respect to a suitable partition of unity $(\varphi_k)_{k\in\Z}$, and equipped with the norm
\[
\|g\|_{\B_{\text{dyad}}} := \sum_{k\in\Z} \|g * \varphi_k\|_\infty < \infty.
\]
The space $\B_{\text{dyad}}$ can be described as follows:
\[
\B_{\text{dyad}} = \left\{g \in \operatorname{BUC}(\R): g = \sum_{k\in\Z} g*\varphi_k, \|g\|_{\B_{\text{dyad}}} :=\sum_{k\in\Z} \|g * \varphi_k\|_\infty < \infty\right\}.
\]
The restriction map and its inverse are both bounded for the norms on the two spaces.  Thus  $\Bq$  is isomorphic to $\B_{\text{dyad}}$.

This relationship between $\B$ and $\B_{\text{dyad}}$ was known and used in \cite{V1}, but the proof there is not very clear.   A precise statement and proof are given in \cite[Proposition 6.2]{BGT}, using the density of $\G$ in $\Bq$ as established in Theorem \ref{Gdense} by means of Arveson's spectral theory for $C_0$-groups.

\subsection*{The $\A$-calculus}
In \cite{ALM}  Arnold and Le Merdy define an $\A$-calculus for negative generators of bounded $C_0$-semigroups on Hilbert spaces, by adapting ideas from \cite{Pel}.  They introduce the space
\[
\A(\R):=  \left\{\sum_{k=1}^\infty g_k*h_k  : g_k \in \operatorname{BUC}(\R), h_k \in H^1(\R),  \sum_{k=1}^\infty \|g_k\|_\infty \|h_k\|_{H^1} < \infty \right\},
\]
with the norm of $f \in\A(\R)$ given by
\[
\|f\|_\A := \inf \left\{\sum_{k=1}^\infty \|g_k\|_\infty \|h_k\|_{H^1} : f = \sum_{k=1}^\infty g_k*h_k, \, g_k \in \operatorname{BUC}(\R), h_k \in H^1(\R)\right\},
\]
where $H^1(\R) = \{f^b: f \in H^1(\C_+)\}$.   Using a quotient of the projective tensor product of the Banach algebras $\operatorname{BUC}(\R)$ (with pointwise multiplication) and $H^1(\R)$ (under convolution), they show that $\A(\R)$ is complete.   Using an idea from \cite{Pel} they prove that $A(\R)$ is an algebra under pointwise multiplication, so it is a Banach algebra which is continuously included in $\operatorname{BUC}(\R) \cap H^\infty(\R)$.    A corresponding Banach algebra $\A(\C_+)$ contained in  $H^\infty(\C_+)$ is obtained from $\A(\R)$  by the Poisson extension process above.    There is a similar construction of spaces $\A_0(\R)$ and $\A_0(\C_+)$, with $\operatorname{BUC}(\R)$ replaced by $C_0(\R)$.    The algebras  $\B$  and $\Bq$ are  densely and continuously included in $\A(\C_+)$ and $\A_0(\C_+)$ respectively.

\begin{thm}
Let $-A$ be the generator of a bounded $C_0$-semigroup on a Hilbert space $X$.  There is a unique bounded homomorphism $\Xi_\A : \A(\C_+) \to L(X)$ which extends the Hille-Phillips calculus and the $\Bes$-calculus for $A$.
\end{thm}


The authors use the half-plane calculus to define $f(A)$ for $ f \in A(\mathbb C_+)$.   
They first obtain a fundamental estimate showing that if 
$f=\mathcal L(\mathcal F^{-1} g \mathcal F^{-1}h)$ with $g \in H^1(\mathbb R)$
and $h, \mathcal F^{-1} h \in L^1(\mathbb R),$ then
$\|f(A)\|\le K_A^2 \|g\|_{L^1}\|f\|_{\infty}$. This depends on approximation
of $H^2(\mathbb R)$-functions by Schwartz functions, factorisation of $H^1$-functions
as products of two $H^2$-functions, and the vector-valued Plancherel theorem (the only place where the Hilbert space structure is used).
By embedding $\A_0(\mathbb R)^*$ into the space of Fourier multipliers
on $H^1(\mathbb R)$, the authors prove that $\L L^1 := \{ \L h: h \in L^1(\R_+)\}$
is dense in $A_0(\mathbb C_+)$.
It follows from this and the fundamental estimate that there is 
a homomorphism from $\A_0(\mathbb C_+)$ into $L(X)$, which coincides
 with the Hille-Phillips calculus on $\L L^1$. Finally the homomorphism extends
to $\A(\mathbb C_+)$ by an approximation procedure similar in spirit to the extension 
of the $\mathcal B$-calculus from $\mathcal{LM}$ to $\mathcal B$ in the second step of the proof of Theorem \ref{besc}. 

\subsection*{White's approach}
White considers negative generators of bounded $C_0$-semigroups on Hilbert spaces in Section 5 of his thesis \cite{White}.   He follows Peller's strategy in \cite{Pel}, but the arguments are more complicated in the continuous case,
and they rely on a duality method and several quite subtle identifications.
Using the injective and projective products of $L^1(\mathbb R_+)$ along with the duals of the products, he derives a continuous-time version of Grothendieck's inequality and obtains several estimates for $f(A)$ involving tensor product norms.  
Passing to Schur multipliers he obtains an estimate
\begin{equation}\label{kernels}
\|f(A)\|\le  K_G K_A^2\sup\left \{\left|\int_{0}^\infty(\mathcal F^{-1} f^b)(t)h(t)\, d\mu(t)\right|:\|h\|_{M({\rm BHK})}\le 1\right \}.
\end{equation}
Here $f \in H^1(\C_+)$, $f^b \in L^1(\R)$ and $\operatorname{supp} (\F^{-1}f^b)$ is a compact subset of $\R_+$, and $h \in L^\infty(\mathbb R_+)$ is a pointwise (Schur) multiplier on a space ${\rm BHK}$ of Hankel type kernels of bounded integral operators on $H^2(\mathbb R),$  $\|\cdot\|_{M({\rm BHK})}$ stands for the multiplier norm, and $K_G$ denotes Grothendieck's constant. 
He constructs an isometric embedding $q$ of BHK into $H^1(\R)^*$, and 
then he shows that $f^b$ belongs to a ``projective convolution
space'', $\mathcal C:=H^1(\mathbb R)*q({\rm BHK})$, which is defined similarly to the space $\A(\mathbb R)$ above. Moreover, $\|f^b\|_{\mathcal C}\le C \|f^b\|_{\mathcal B_{\rm dyad}}$, where $C$ is an absolute constant.   
From these identifications, the estimate \eqref{kernels} and the isomorphism between $\mathcal B_{\rm dyad}$ and $\mathcal B_0$, it follows that there is an absolute constant $C_1$ such that
\[
\|f(A)\|\le C_1 K_A^2 \|f\|_{\mathcal B_0}
\]
if $f \in \G$ (as in \eqref{Gest}) and $f^b$ is a Schwartz function on $\R$.
This estimate does not appear explicitly in the thesis, but it can be seen from the arguments leading to \cite[Theorem 5.5.12]{White}.  

\subsection*{Cayley transform question}

The various estimates in \cite{ALM} and \cite[Section 5]{White} involve new norms on $\B$ which are dominated by the $\B$-norm, and which provide upper bounds for rational functions of $A$.  This is relevant to the Cayley transform question on Hilbert spaces in the following way. Here we will  assume that the norm is the $\A$-norm, but the same comments apply to other norms of this type.   

If $-A$ generates a bounded $C_0$-semigroup, then it follows from the boundedness of the Hille-Phillips calculus with respect to the $\A$-norm that
\[
\|V(A)^n\| \le  C  \|V^n\|_\A, \qquad n\in\N,
\]
where $C$ is a constant.   Thus it would be interesting to determine the asymptotic behaviour of the norms $\|V^n\|_\A$ for large $n$.  As seen in Application 4 of Section \ref{B-calc}, they cannot grow faster than logarithmically, and they are bounded below by $1/C$.  If the norms are bounded above by a constant, then the Cayley transform question for Hilbert spaces has a positive answer.   On the other hand, no new information is obtained if the norms grow logarithmically.   There is also the possibility that the norms grow at a rate which is slower than logarithmic. 

Unfortunately, the $\A$-norms of simple functions such as $V^n$ are difficult to estimate, because of the indirect way that the $\A$-norm is defined.   Moreover we do not know of a simple reproducing formula for the algebra $\A$ of the style that we have for $\B$, $\D$ and $\H$.   Thus the asymptotic behaviour of $\|V^n\|_\A$ for large $n$ remains unknown to us. This also applies to the other norms of this type, but not to the $\B$-norm for which the formula \eqref{Bnorm} led to a sharp estimate.





\end{document}